\documentclass[letter,11pt, reqno]{amsart}

\usepackage{natbib}
 \bibpunct[, ]{(}{)}{,}{a}{}{,}%

 \usepackage[hidelinks,colorlinks = true,linkcolor={blue!80!black},citecolor={blue!50!black},urlcolor={blue!80!black}]{hyperref}
\usepackage{amsthm,amsmath,amsfonts}
\usepackage{amssymb}
\usepackage{url}
\usepackage{lineno}
\usepackage[linesnumbered,lined,boxed,commentsnumbered]{algorithm2e}

\usepackage{subcaption}
\usepackage{adjustbox}
\usepackage{tabularx,booktabs,multirow}
\usepackage{booktabs}
\usepackage{multicol}
\usepackage{soul}
\usepackage{bm}
 
\usepackage{tikz}
\usetikzlibrary{positioning,calc}
\usetikzlibrary{arrows.meta, decorations.pathreplacing}

\usepackage{fontawesome5}
\usepackage[T1]{fontenc}
\usepackage{xcolor}

\definecolor{spolddraw}{RGB}{180,210,255}
\definecolor{spoldfill}{RGB}{220,235,255}

\definecolor{spnewdraw}{RGB}{40,90,200}
\definecolor{spnewfill}{RGB}{74,144,246}

\definecolor{smolddraw}{RGB}{245,190,110}
\definecolor{smoldfill}{RGB}{255,228,181}

\definecolor{smnewdraw}{RGB}{230,120,20}
\definecolor{smnewfill}{RGB}{255,170,70}

\definecolor{reactionolddraw}{RGB}{120,120,120}
\definecolor{reactionoldfill}{RGB}{220,220,220}

\definecolor{reactionnewdraw}{RGB}{200,60,60}
\definecolor{reactionnewfill}{RGB}{255,190,190}

\def\R{\mathbb{R}}

\def\A{\mathcal{A}}
\def\N{\mathcal{N}}

\newtheorem{definition}{Definition}[section]

\usepackage[colorinlistoftodos]{todonotes}
\usepackage{natbib}
\definecolor{cerulean}{rgb}{0.0, 0.40, 0.60}

\definecolor{darkgreen}{RGB}{34, 139, 34}

\usepackage[english]{babel}

\usepackage[letterpaper,top=2cm,bottom=2cm,left=3cm,right=3cm,marginparwidth=1.75cm]{geometry}

\usepackage{amsmath}
\usepackage{graphicx}

\let\origmaketitle\maketitle
\def\maketitle{
	\begingroup
	\def\uppercasenonmath##1{} % this disables uppercasing title
	\let\MakeUppercase\relax % this disables uppercasing authors
	\origmaketitle
	\endgroup
}

\title{\Large Constructing Nested Self-Amplifying Multiperiod Hypergraphs through Mathematical Optimization}

\author[V. Blanco, R. G\'azquez, \MakeLowercase{and} J.F. Oca\~na-Rivas]{
{\large V\'ictor Blanco$^{\dagger,\ddagger}$, Ricardo G\'azquez$^{\dagger,\ddagger}$, Juan Francisco Oca\~na-Rivas$^{\dagger}$}\medskip\\
$^\dagger$Institute of Mathematics (IMAG), Universidad de Granada\\
$^\ddagger$Dpt. Quantitative Methods for Economics \& Business, Universidad de Granada\\
\texttt{vblanco@ugr.es}, \texttt{rgazquez@ugr.es}, \texttt{jfocana@correo.ugr.es}\\
}


\begin{document}

\begin{abstract}
This paper proposes an optimization-based framework for the analysis of multiperiod directed multihypergraphs aimed at identifying self-amplifying structures that sustain endogenous growth in complex systems. The approach captures the progressive and nested activation of nodes and hyperarcs, providing a dynamic representation of evolving production and reaction networks. We formulate the problem as a mixed integer optimization model. First, we introduce a tractable linear formulation that captures structural amplification. We then extend this model to a mixed integer nonlinear setting that incorporates a \emph{synergistic flow law} that generalizes mass-action kinetics in Chemical Reaction Networks and that accounts for interaction effects. This nonlinear formulation is handled through logarithmic transformations and piecewise-linear outer approximations. The framework unifies combinatorial structure selection and flow dynamics, bridging Mathematical Optimization with applications in Economics and Chemistry, including autocatalytic systems related to the \emph{Origin of Life}. Computational experiments on synthetic instances demonstrate scalability, while an input--output case study illustrates the ability of the model to identify growth-enabling sectors, interdependencies, and structural bottlenecks across different periods, providing actionable insights for the analysis and management of complex systems.
\end{abstract}

\maketitle

\section{Introduction}

Understanding how complex systems grow, evolve, and self-amplify over time is a central question in economics, operations research, and management science. In many settings, such as production economies, supply chains, innovation systems, and industrial processes, growth is not merely the result of exogenous inputs, but rather emerges from \emph{internal feedback mechanisms} in which outputs are reinvested as inputs to sustain and expand activity~\citep{aghion2008economics,huang2014ecological}. Classical models of economic growth, dating back to \cite{neumann1945model}, formalize this idea through input--output structures and technological systems capable of generating endogenous expansion. More recent work has further highlighted how decentralized interactions and adaptive dynamics can lead to sustained growth, while simultaneously producing structural imbalances such as inequality \citep{branzei2018universal,acemoglu2012network,acemoglu2017networks,baqaee2019granular,liu2024introduction,nagurney2013supply,wongthatsanekorn2010multi}, particularly in networked and platform-based environments where local interactions generate aggregate amplification effects.

A key feature underlying these phenomena is the presence of \emph{self-reinforcing production structures}, where a subset of activities or technologies collectively generates more resources than it consumes. These structures can be interpreted as economic analogues of amplification mechanisms: they identify combinations of processes that are capable of sustaining long-term growth. From an analytical perspective, such systems can be naturally modeled using network-based representations, where nodes correspond to goods, resources, or activities, and interactions capture production or transformation technologies relationships. However, many real-world processes involve \emph{multi-input, multi-output interactions}, which cannot be adequately represented by pairwise networks, thus requiring more general modeling frameworks such as hypergraphs or reaction networks~\citep{baqaee2019granular,benati2025modularity,feinberg2019foundations}.

Hypergraphs provide a powerful modeling framework for such systems, allowing each interaction to involve multiple inputs and outputs simultaneously. This is particularly relevant in economic and industrial contexts, where production technologies combine several resources to generate multiple products, and where complex interdependencies arise across sectors. Recent advances have shown that self-reinforcing structures in hypergraphs can be rigorously characterized through optimization-based approaches, leading to quantitative measures such as the maximal amplification factor and efficient methods to identify growth-enabling subsystems \citep{blanco2025identifying,gagrani2023geometry}. These developments connect naturally with classical economic growth theory, as they extend von Neumann-type expansion models to general multi-input production networks. 

While these approaches provide powerful tools to detect \emph{whether} a system can sustain growth, they are inherently static: they identify feasible growth structures but do not explain \emph{how such structures emerge over time}. While these approaches provide powerful tools to detect \emph{whether} a system can sustain growth, they are inherently static: they identify feasible growth structures but do not explain \emph{how such structures emerge over time}, a limitation that has been widely recognized in the literature on path dependence and dynamic network formation \citep{arthur1989competing}. In real economic and operational systems, growth is a dynamic and path-dependent process. New technologies, production activities, or supply chain links are gradually introduced, often depending on the prior availability of resources or intermediate goods. Similarly, in decentralized systems, local decisions and interactions shape the trajectory of growth, leading to heterogeneous expansion patterns and structural asymmetries.

Motivated by these considerations, this paper introduces the \emph{Growth and Emergence Model} (GEM), a novel optimization framework for modeling the \emph{temporal formation of self-amplifying structures} in directed multihypergraphs. Rather than identifying a single static production structure, GEM constructs a sequence of nested subhypergraphs that represent the progressive activation of activities and interactions across multiple periods. This provides a dynamic and interpretable representation of how growth unfolds, capturing both structural feasibility and temporal evolution. 

A key advantage of GEM is its ability to generate \emph{explainable growth trajectories}. By explicitly modeling the sequence of activation decisions, the framework reveals which components drive expansion, how interdependencies shape growth, and how resources propagate through the system over time. This is particularly valuable in economic and managerial contexts, where understanding the drivers of growth and the role of different activities is essential for decision-making, planning, and policy design.

Building on this framework, we propose two versions of the model. First, we introduce \emph{GEM-E}, an equitable formulation that incorporates an ordered utility objective to promote balanced activation across periods. This objective belongs to the class of \emph{ordered aggregation operators}, most notably the ordered weighted averaging (OWA) operator introduced by \citet{Yager1988}. OWA operators expand modeling capabilities by capturing preferences related to robustness, fairness, and distributional sensitivity. Specifically, the OWA operator reorders a vector of values and computes a weighted average based on their rank, thereby generalizing classical criteria such as the minimum, maximum, median, quantiles, and mean.

This flexibility has motivated a broad literature in operations research and optimization, particularly in ordered location and network design problems~\citep[e.g.,][]{blanco2023fairness,cherkesly2025ranking,Marin2020,Puerto2019,Ljubic2024,ogryczak2016ordered}, as well as applications in insurance~\citep{Blanco2018BMS}, portfolio optimization~\citep{BENATI2025103274}, network design~\citep{Puerto2015OrderedMH,disc:blanco2019ordered}, or robust regression~\citep[e.g.,][]{Yager2009}, among many others. In the context of GEM, the OWA-based objective enables the construction of growth trajectories that balance efficiency and equity across periods, avoiding solutions in which activation is overly concentrated in a small subset of stages and leading instead to more stable and interpretable expansion patterns.

Second, we develop \emph{GEM-D}, a dynamics-aware version of GEM-E that incorporates a \emph{Synergistic Flow Law} together with discrete-time state dynamics. This formulation captures how interactions among activities influence future growth by linking activation decisions across periods through endogenous flow relationships. The \emph{Synergistic Flow Law} extends classical linear flow models by incorporating complementarities across inputs, whereby the joint availability of resources amplifies production outcomes. This mechanism is analogous to mass-action kinetics in chemical systems and to multiplicative interaction models widely used in economics, such as Cobb--Douglas production functions and matching models, where output depends on the simultaneous presence of multiple inputs or agents. As a result, GEM-D bridges optimization-based modeling with dynamic system behavior, providing a richer and more realistic representation of growth processes driven by interaction and reinforcement mechanisms. 

Although the proposed framework is novel in the context of growth and expansion in hypergraphs, the modeling of multiperiod decisions has been extensively studied in a wide range of combinatorial optimization problems with social relevance. In particular, dynamic and multiperiod formulations have been successfully applied in areas such as network design, facility location, and resource allocation, where decisions must account for temporal evolution and inter-period dependencies \citep[see, e.g.][]{albareda2014dynamic,BELTRANROYO201626,marin2018multi,nickel2020multi}.

Although our primary motivation arises from economic and operational systems, the proposed framework is sufficiently general to encompass other domains characterized by self-amplifying interactions. In particular, chemical reaction networks and biological systems provide natural applications, where the emergence of autocatalytic pathways reflects similar amplification mechanisms governed by network structure, kinetics, and thermodynamic constraints \citep{feinberg2019foundations,deshpande2014autocatalysis,blokhuis_universal_2020}. 

A substantial body of work has shown how autocatalytic and self-replicating structures can spontaneously arise in reaction systems and play a fundamental role in the organization of biochemical processes and, more fundamentally, in theories on the \emph{Origin of Life}, where the emergence of self-sustaining chemical networks is considered a key step in the transition from chemistry to biology \citep{liu2018mathematical,peng2020ecological,hordijk2018autocatalytic,andersendefining,peng2022wilhelm}. These studies, together with recent optimization-based approaches for identifying amplification structures \citep{gagrani2025thermodynamic,gagrani2026topological}, highlight the deep connections between network topology, interaction rules, and emergent growth dynamics. 

This cross-domain perspective underscores the broad applicability of the proposed approach and reinforces the fundamental role of network structure in governing growth processes across economic, biological, and technological systems. 

From a methodological standpoint, the resulting models lead to mixed-integer nonlinear optimization formulations that integrate hypergraph structure, flow variables, and temporal constraints. To ensure tractability, we derive equivalent reformulations and propose piecewise-linear approximations that can be efficiently solved with modern optimization tools.  This enables the application of GEM to realistic instances arising in production systems, supply chain design, and networked economic environments.

The contributions of this paper are fourfold. First, we introduce GEM, a novel framework for modeling the temporal emergence of self-reinforcing production structures. Second, we develop GEM-E, which incorporates equity considerations into the growth process. Third, we propose GEM-D, which integrates dynamic interactions and system evolution. Fourth, we provide tractable formulations and computational strategies that enable practical implementation in large-scale settings. The proposed methodologies are first validated on synthetic (yet realistic) hypergraphs generated using a random chemical reaction network generator (\texttt{SMGen}), and subsequently applied to a real-world case study based on Bureau of Economic Analysis (BEA) data, from which managerial insights on the dynamics of economic activities and technologies are derived.

The rest of the paper is organized as follows. 
Section~\ref{sec:1} introduces the notation and core concepts on directed multihypergraphs, including incidence structures, hyperflows, and the definition of self-reinforcing substructures that underpin the proposed framework. 
Section~\ref{sec:2} presents the GEM-E formulation, detailing the ordered utility objective, the equity-driven modeling features, and the constraints ensuring feasibility and nested growth across periods. 
Section~\ref{sec:3} develops the GEM-D model, incorporating the Synergistic Flow Law and discrete-time dynamics, and provides tractable reformulations and approximation strategies for its efficient solution. 
Section~\ref{sec:4} reports numerical experiments on synthetic instances inspired by chemical reaction networks. Section \ref{sec:5} is devoted to illustrate the results obtained in our economic case study illustrating the applicability of the proposed models in production systems using the BEA dataset. 
Finally, Section~\ref{sec:6} concludes with a summary of the main findings, discusses managerial and methodological implications, and outlines directions for future research.
\section{Preliminaries}\label{sec:1}

In this section, we introduce the notation and fundamental concepts underlying the \emph{Growth and Emergence Model} (GEM). These elements provide the mathematical foundation for representing systems in which growth arises from the interaction of multiple components over time.

Our primary combinatorial object is the \emph{directed multihypergraph}, which generalizes classical network models by allowing interactions among multiple entities simultaneously. This representation is particularly suitable for modeling production systems, technological processes, and networked economies, where activities combine several inputs to generate multiple outputs. More generally, directed multihypergraphs provide a flexible framework for describing higher order dependencies in complex systems~\citep{klamt2009hypergraphs,battiston2021physics,bick2023higher}.

\begin{definition}
    A \emph{directed multihypergraph} is a tuple $\mathcal{H} = (\mathcal{N}, \A)$, where $\mathcal{N}$ is the set of nodes and $\A$ is the set of hyperarcs. Each hyperarc $a \in \A$ is represented as an ordered pair $a = (S_a, T_a)$, where $\emptyset \neq S_a, T_a \subseteq \mathcal{N}$ denote the source and target multisets, respectively. Multiplicities within $S_a$ and $T_a$ allow nodes to appear multiple times, modeling intensity, stoichiometry, or input-output coefficients.
\end{definition}

In what follows, we use the term \emph{hypergraph} to refer to the more general case of a \emph{multihypergraph}. The distinction lies in node multiplicities: in a hypergraph, each node appears with multiplicity one in every hyperedge, whereas in a multihypergraph arbitrary multiplicities are allowed. By a slight abuse of terminology, we use the term hypergraph to denote this general setting throughout the paper.

To illustrate this object, consider a simple economy with three goods: 
raw material ($R$), intermediate product ($I$), and final good ($F$). 
We model the production technologies as a directed hypergraph 
$\mathcal{H} = (\mathcal{N}, \A)$, where 
$\mathcal{N} = \{R, I, F\}$ represents the set of available goods in the economy and $\A = \{a_1, a_2, a_3\}$ are the different technologies. Specifically:
    \begin{align*}
        a_1 &= (\{R, R\}, \{I\}) 
        && \text{(two units of raw material produce one intermediate product)}, \\
        a_2 &= (\{I\}, \{F, F\}) 
        && \text{(one intermediate produces two final goods)}, \\
        a_3 &= (\{F\}, \{R\}) 
        && \text{(recycling: one final good yields one unit of raw material)}.
    \end{align*}

\noindent
This hypergraph models a simple production cycle: raw materials are transformed into 
intermediate goods, which are then converted into final products, while part of the 
output can be recycled back into raw materials. Such structures naturally arise in 
input--output economic models and generalize classical production systems by allowing 
multi-input and multi-output technologies.

For each node $v \in \mathcal{N}$ and hyperarc $a \in \A$, define
$$
\mathbb{S}_{va} = \text{multiplicity of } v \text{ in } S_a, 
\qquad 
\mathbb{T}_{va} = \text{multiplicity of } v \text{ in } T_a.
$$
The source and target incidence matrices are
$$
\mathbb{S} = (\mathbb{S}_{va}), \qquad \mathbb{T} = (\mathbb{T}_{va}),
$$
and the net incidence matrix is $\mathbb{Q} = \mathbb{T} - \mathbb{S}$, where the entry $\mathbb{Q}_{va}$ represents the net contribution of node $v$ in hyperarc $a$. Negative values correspond to consumption, positive values to production.

In the example above, the incidence matrices are:
$$
\mathbb{S} = \begin{pmatrix}
    2 &0 &0\\
    0 &1&0\\
    0&0&1
\end{pmatrix}, \quad \mathbb{T} =\begin{pmatrix}
    0 &0&1\\
    1 &0&0\\
    0&2&0
\end{pmatrix}.
$$

A hyperflow is a vector $\mathbf{f} \in \mathbb{R}^{|\A|}_{>0}$ assigning an activity level to each hyperarc. The net balance is given by $\mathbf{x} = \mathbb{Q}\mathbf{f} \in \mathbb{R}^{|\mathcal{N}|}$, 
where negative values in $\mathbf{x}$ represent \emph{depletion} and positive values represent \emph{amplification} of the nodes.

If, in the example above, we consider a flow vector $\mathbf{f} = (1, 3, 2)$, 
indicating an intensity of $1$ unit for activity $a_1$, $3$ units for $a_2$, 
and $2$ units for $a_3$, the resulting net production of the system is given by
$$
\mathbf{x} = \mathbb{Q}\mathbf{f} = 
\begin{pmatrix}
    0 \\ 
    -2 \\ 
    4
\end{pmatrix}.
$$
This implies that, under the intensity vector $\mathbf{f}$, there is no net 
production of good $R$, a net depletion of two units of good $I$, and a net 
production of four units of the final good $F$.

A subhypergraph $\mathcal{H}' = (\mathcal{N}', \A')$ satisfies $\A' \subseteq \A$ and $S_a, T_a \subseteq \mathcal{N}'$ for all $a \in \A'$. Given $\mathcal{M} \subseteq \mathcal{N}'$, the restricted subhypergraph $\mathcal{H}'|_{\mathcal{M}}$ contains all hyperarcs whose source and target multisets are fully contained in $\mathcal{M}$.

\begin{definition}[Self-Amplifying Subhypergraph]\label{def:selfamplifying}
Let $\mathcal{H} = (\mathcal{N}, \A)$ be a directed hypergraph and let $\mathcal{M} \subseteq \mathcal{N}$. 
Let $\mathcal{H}'=(\mathcal{N}', \mathcal{A}') \subseteq \mathcal{H}$ a subhypergraph of $\mathcal{H}$ with $\mathcal{M}\subseteq \mathcal{N}'$.

We say that the tuple $(\mathcal{H'};\mathcal{M})$ defines a \textbf{self-amplifying subhypergraph} if the following conditions hold:

\begin{enumerate}
\item \textbf{Self-sufficiency:}
\begin{align*}
&\forall a \in \A',\ \exists v,v'\in\mathcal{M}:\ \mathbb{T}_{va}>0,\ \mathbb{S}_{v'a}>0,\\
&\forall v \in \mathcal{M},\ \exists a,a'\in\A':\ \mathbb{T}_{va}>0,\ \mathbb{S}_{va'}>0.
\end{align*}

\item \textbf{Net positive realizability:}
$$
\exists \mathbf{f} \in \mathbb{R}_{>0}^{|\A'|} \text{ such that } 
\sum_{a\in\A'} \mathbb{Q}_{va} f_a > 0, \quad \forall v \in \mathcal{M}.
$$
\end{enumerate}
\end{definition}
Nodes in $\mathcal{M}$ are called self-amplifying nodes. A self-amplifying subhypergraph is minimal if it contains no proper self-amplifying substructure.

In the example above, the selected hypergraph, together with the flow $\mathbf{f}$, does not constitute a fully self-amplifying structure, since not all components are simultaneously expanded (node $I$ exhibits a depletion of $2$ units, while node $R$ has zero net production). Therefore, the net positive realizability condition is not satisfied. In fact, it is not possible to find a flow that satisfy the net realizability for $\mathcal{M}=\{R,I,F\}$ since the system:
\begin{align*}
    \mathbb{Q} \mathbf{f} > 0 \Rightarrow \begin{cases}
        -2 f_1 + f_3 >0\\
        f_1 - f_2 >0\\
        2f_2 - f_3 >0
    \end{cases} \Rightarrow f_3 > 2f_1 > 2f_2 > f_3
\end{align*}
is clearly infeasible.

However, the structure defined by $\mathcal{H}' = \mathcal{H}$, $\mathcal{M} = \{R,I\}$, and $\tilde{\mathbf{f}} = (2,1,5)$ does satisfy the self-amplifying condition. In particular, both $R$ and $I$ fulfill the self-sufficiency property, as they appear as both input and output nodes within the hypergraph, and every hyperarc in $\mathcal{A}$ involves at least one of these nodes in its input and output sets.  Moreover, the net positive realizability condition holds, since
$$
\sum_{a \in \mathcal{A}} \mathbb{Q}_{Ra} \tilde{f}_a = 1 > 0 
\quad \text{and} \quad 
\sum_{a \in \mathcal{A}} \mathbb{Q}_{Ia} \tilde{f}_a = 1 > 0.
$$
\noindent Hence, under the hyperflow $\tilde{\mathbf{f}}$, both $R$ and $I$ are strictly positively produced.

This example highlights that self-amplification is not a property of the hypergraph alone, but rather of the combination of a subset of nodes and a feasible hyperflow. In particular, suitable restrictions of the system may reveal internally self-sustaining structures even when the full system does not exhibit global amplification. However, this observation also reveals a significant computational challenge: identifying such structures in large hypergraphs requires jointly selecting a subset of nodes and determining a compatible hyperflow that satisfies the self-amplification conditions. The number of possible subsets grows exponentially with the size of the system, and for each candidate subset, feasibility must be verified through a system of inequalities involving the hyperflow variables. As a result, the detection of self-amplifying subhypergraphs naturally leads to a combinatorial optimization problem of high complexity, motivating the need for systematic optimization-based approaches capable of efficiently exploring the space of feasible substructures, as those already proposed by  \citet{peng_hierarchical_2022}, \citet{gagrani2023geometry}, or \citet{blanco2025identifying}.

Given a time horizon $T\in \mathbb{Z}_{>0}$, we define the corresponding index set as $\mathcal{T} := \{1,\ldots, T\}$. The framework then constructs sequences of nested self-amplifying chain of subhypergraphs over time.
\begin{definition}
A sequence $\{(\mathcal{H}^t; \mathcal{M}^t)\}_{t\in \mathcal{T}}$ with $\mathcal{H}^t=(\mathcal{N}^t,\A^t)$ is a nested self-amplifying  chain of subhypergraphs if
$$
\mathcal{H}^1 \subseteq \cdots \subseteq \mathcal{H}^T,
$$
\noindent and the subset $\mathcal{M}^t$ that is self-amplifying. 
\end{definition}

In our first framework, denoted as \emph{GEM-E} (Section \ref{sec:2}), we focus on constructing such chains through mathematical optimization models that ensure both the nested structure of the subhypergraphs and the self-amplifying condition at every period. The second framework builds upon this approach by incorporating additional interactions between the flows that ensure the net positive realizability of these subhypergraphs. These interactions are governed by the so-called \emph{Synergistic Flow Law}, which is defined as follows.

\begin{definition}[Synergistic Flow Law]\label{def:synergistic}
Let $\mathcal{H} = (\mathcal{N}, \A)$ be a directed hypergraph. Given $\boldsymbol{\kappa} = (\kappa_a)_{a \in \A} \in \mathbb{R}_{>0}^{|\A|}$ (synergy coefficients), we say that a hyperflow $\mathbf{f} \in \mathbb{R}_{> 0}^{|\A|}$ 
satisfies the \textbf{Synergistic Flow Law} if there exists a vector of node activity levels $\mathbf{x} = (x_v)_{v \in \mathcal{N}} \in \mathbb{R}^{|\mathcal{N}|}$  such that:
\begin{equation}\label{ctr:12_0}
    f_a = \kappa_a \prod_{v \in \mathcal{N}} x_v^{\mathbb{S}_{va}}, 
\quad \forall a \in \A.
\end{equation}
\end{definition}
This law generalizes multiplicative interaction models such as Cobb-Douglas production functions, where output depends on the joint availability of inputs.

Note that a self-amplifying subhypergraph guarantees the existence of a hyperflow $\mathbf{f}$ satisfying the net positive realizability condition for the self-amplifying nodes. When such a flow also complies with the Synergistic Flow Law introduced above, we refer to the subhypergraph as \textit{synergistic self-amplifying}.  This notion naturally extends to nested chains of subhypergraphs by requiring the Synergistic Flow Law to hold across all periods. The construction of such structures constitutes the objective of our second framework, denoted as \emph{GEM-D} (Section \ref{sec:3}).

To illustrate the two frameworks, we consider a small synthetic hypergraph with seven nodes and six hyperarcs with the following  incidence matrices and synergy coefficients:
$$
\mathbb{S} =
\begin{pmatrix}
1 & 0 & 1 & 0 & 0 & 0 \\
0 & 1 & 0 & 0 & 0 & 0 \\
0 & 0 & 0 & 1 & 1 & 1 \\
0 & 0 & 0 & 0 & 0 & 0 \\
0 & 0 & 0 & 0 & 0 & 1 \\
0 & 1 & 1 & 0 & 0 & 0 \\
0 & 0 & 0 & 0 & 1 & 0
\end{pmatrix},
\qquad
\mathbb{T} =
\begin{pmatrix}
0 & 0 & 0 & 1 & 0 & 0 \\
1 & 0 & 0 & 0 & 1 & 1 \\
0 & 0 & 1 & 0 & 0 & 0 \\
1 & 0 & 0 & 0 & 0 & 0 \\
0 & 0 & 0 & 1 & 0 & 0 \\
0 & 1 & 1 & 0 & 0 & 1 \\
0 & 0 & 0 & 0 & 0 & 0
\end{pmatrix},
\qquad
\kappa = \begin{pmatrix}
9\\
4.4\\
6.7\\
4.0\\
1.5\\
6.8
\end{pmatrix}.
$$

\begin{figure}[htbp]
\centering
\begin{subfigure}{0.47\textwidth}
\centering
\resizebox{\linewidth}{!}{\fbox{\begin{tikzpicture}[>=Latex, line cap=round, line join=round, x=2cm, y=1.2cm]

% ---- estilos ----
\tikzset{
  species/.style={draw, circle, minimum size=6mm, inner sep=1pt, font=\scriptsize},
  species-inactive/.style={species, fill=white},
  species-oldactive/.style={species, fill=orange!20},
  species-newactive/.style={species, fill=orange!70},
  species-out-oldactive/.style={species, fill=blue!20},
  species-out-newactive/.style={species, fill=blue!70},
  reaction-active/.style={draw=red, rectangle, minimum size=5mm, inner sep=2pt, font=\scriptsize, line width=0.8pt},
  reaction-inactive/.style={draw=black, rectangle, minimum size=5mm, inner sep=2pt, font=\scriptsize, fill=white},
  species-hidden/.style={species, draw=none, fill=none},
  reaction-hidden/.style={rectangle, minimum size=5mm, inner sep=2pt, draw=none, fill=none},
  edge-in-active/.style={->, draw=orange!80, line width=0.8pt},
  edge-in-inactive/.style={->, draw=orange!30, line width=0.6pt},
  edge-out-active/.style={->, draw=blue!90, line width=0.8pt},
  edge-out-inactive/.style={->, draw=blue!30, line width=0.6pt}
}

% --------------------
% NODOS IZQUIERDA
% --------------------
\node[species-newactive] (S0L) at (0,6) {S1};
\node[species-newactive] (S1L) at (0,5) {S2};
\node[species-newactive] (S2L) at (0,4) {S3};
\node[species-hidden] (S3L) at (0,3) {};
\node[species-newactive] (S4L) at (0,2) {S5};
\node[species-newactive] (S5L) at (0,1) {S6};
\node[species-hidden] (S6L) at (0,0) {};

% --------------------
% REACCIONES
% --------------------
\node[reaction-hidden] (R0) at (2,5.5) {};
\node[reaction-active] (R1) at (2,4.5) {R2};
\node[reaction-active] (R2) at (2,3.5) {R3};
\node[reaction-active] (R3) at (2,2.5) {R4};
\node[reaction-active] (R4) at (2,1.5) {R5};
\node[reaction-active] (R5) at (2,0.5) {R6};

% --------------------
% NODOS DERECHA
% --------------------
\node[species-out-newactive] (S0R) at (4,6) {S1};
\node[species-out-newactive] (S1R) at (4,5) {S2};
\node[species-out-newactive] (S2R) at (4,4) {S3};
\node[species-hidden] (S3R) at (4,3) {};
\node[species-out-newactive] (S4R) at (4,2) {S5};
\node[species-out-newactive] (S5R) at (4,1) {S6};
\node[species-inactive] (S6R) at (4,0) {S7};

% ---- aristas izquierda -> reacción (misma estructura) ----
\draw[edge-in-active] (S0L) -- (R3);
\draw[edge-in-active] (S1L) -- (R4);
\draw[edge-in-active] (S1L) -- (R5);
\draw[edge-in-active] (S2L) -- (R2);
\draw[edge-in-active] (S4L) -- (R3);
\draw[edge-in-active] (S5L) -- (R1);
\draw[edge-in-active] (S5L) -- (R2);
\draw[edge-in-active] (S5L) -- (R5);

% ---- aristas reacción -> derecha (misma estructura) ----
\draw[edge-out-active] (R1) -- (S1R);
\draw[edge-out-active] (R1) -- (S5R);

\draw[edge-out-active] (R2) -- (S5R);
\draw[edge-out-active] (R2) -- (S0R);

\draw[edge-out-active] (R3) -- (S2R);

\draw[edge-out-active] (R4) -- (S2R);
\draw[edge-out-active] (R4) -- (S6R);

\draw[edge-out-active] (R5) -- (S2R);
\draw[edge-out-active] (R5) -- (S4R);

\end{tikzpicture}}}
\caption{Self-amplifyng subhypergraph.}\label{f:toy_example_t0_MAK0}
\end{subfigure}
\hfill
\begin{subfigure}{0.47\textwidth}\centering
\resizebox{\linewidth}{!}{\fbox{\begin{tikzpicture}[>=Latex, line cap=round, line join=round, x=2cm, y=1.2cm]

% ---- estilos ----
\tikzset{
  species/.style={draw, circle, minimum size=6mm, inner sep=1pt, font=\scriptsize},
  species-inactive/.style={species, fill=white},
  species-oldactive/.style={species, fill=orange!20},
  species-newactive/.style={species, fill=orange!70},
  species-out-oldactive/.style={species, fill=blue!20},
  species-out-newactive/.style={species, fill=blue!70},
  reaction-active/.style={draw=red, rectangle, minimum size=5mm, inner sep=2pt, font=\scriptsize, line width=0.8pt},
  reaction-inactive/.style={draw=black, rectangle, minimum size=5mm, inner sep=2pt, font=\scriptsize, fill=white},
  species-hidden/.style={species, draw=none, fill=none},
  reaction-hidden/.style={rectangle, minimum size=5mm, inner sep=2pt, draw=none, fill=none},
  edge-in-active/.style={->, draw=orange!80, line width=0.8pt},
  edge-in-inactive/.style={->, draw=orange!30, line width=0.6pt},
  edge-out-active/.style={->, draw=blue!90, line width=0.8pt},
  edge-out-inactive/.style={->, draw=blue!30, line width=0.6pt}
}

% --------------------
% NODOS IZQUIERDA
% --------------------
\node[species-hidden] (S0L) at (0,6) {};
\node[species-newactive] (S1L) at (0,5) {S2};
\node[species-newactive] (S2L) at (0,4) {S3};
\node[species-hidden] (S3L) at (0,3) {};
\node[species-hidden] (S4L) at (0,2) {};
\node[species-newactive] (S5L) at (0,1) {S6};
\node[species-hidden] (S6L) at (0,0) {};

% --------------------
% REACCIONES
% --------------------
\node[reaction-hidden] (R0) at (2,5.5) {};
\node[reaction-active] (R1) at (2,4.5) {R2};
\node[reaction-active] (R2) at (2,3.5) {R3};
\node[reaction-hidden] (R3) at (2,2.5) {};
\node[reaction-active] (R4) at (2,1.5) {R5};
\node[reaction-active] (R5) at (2,0.5) {R6};

% --------------------
% NODOS DERECHA
% --------------------
\node[species-inactive] (S0R) at (4,6) {S1};
\node[species-out-newactive] (S1R) at (4,5) {S2};
\node[species-out-newactive] (S2R) at (4,4) {S3};
\node[species-hidden] (S3R) at (4,3) {};
\node[species-inactive] (S4R) at (4,2) {S5};
\node[species-out-newactive] (S5R) at (4,1) {S6};
\node[species-inactive] (S6R) at (4,0) {S7};

% =========================================================
% ARISTAS izquierda -> reacción
% =========================================================
\draw[edge-in-active] (S2L) -- (R2);
\draw[edge-in-active] (S5L) -- (R1);
\draw[edge-in-active] (S5L) -- (R2);
\draw[edge-in-active] (S5L) -- (R5);
\draw[edge-in-active] (S1L) -- (R5);
\draw[edge-in-active] (S1L) -- (R4);

% =========================================================
% ARISTAS reacción -> derecha
% =========================================================
\draw[edge-out-active] (R1) -- (S1R);
\draw[edge-out-active] (R1) -- (S5R);

\draw[edge-out-active] (R2) -- (S5R);
\draw[edge-out-active] (R2) -- (S0R);

\draw[edge-out-active] (R4) -- (S6R);
\draw[edge-out-active] (R4) -- (S2R);

\draw[edge-out-active] (R5) -- (S2R);
\draw[edge-out-active] (R5) -- (S4R);

\end{tikzpicture}}}
\caption{Synergistic self-amplifyng subhypergraph.}
\end{subfigure}
\caption{Single period structure.}
\label{f:toy_example_t0}
\end{figure}
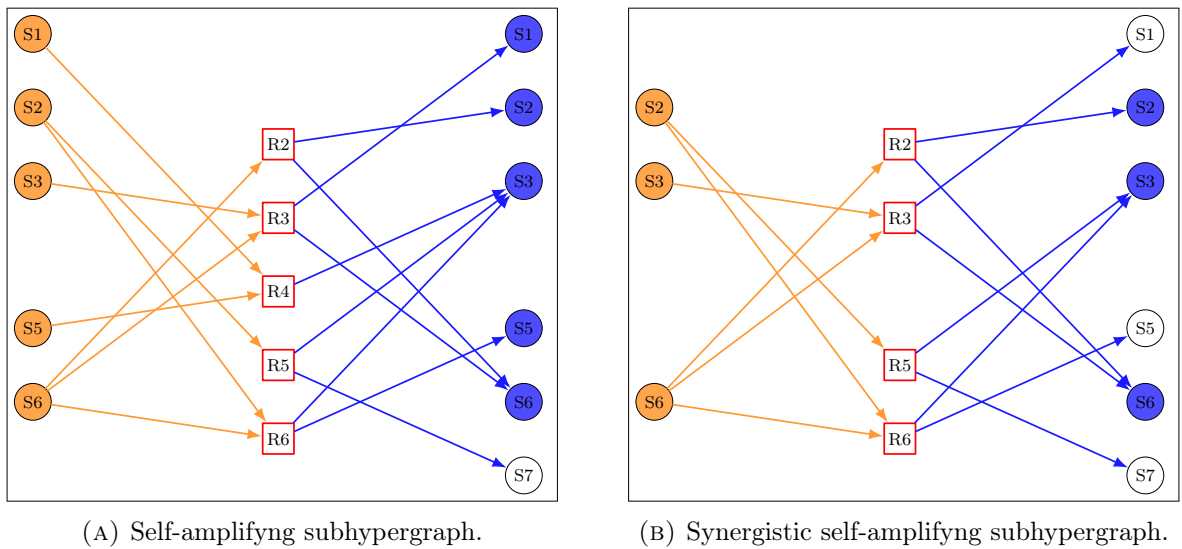

Figure~\ref{f:toy_example_t0} depicts the resulting hypergraph in a single period under the two frameworks (self-amplifying on the left and synergistic self-amplifying on the right). We represent the hypergraph using a tripartite structure in which the nodes (circles) are duplicated (left for inputs and right for outputs), while hyperarcs are represented as squares in the middle. An arrow is drawn from a left node $v$ to a hyperarc $a$ if $v$ is an input node of $a$ (i.e., if $\mathbb{S}_{va}>0$), and an arrow is drawn from the hyperarc to a right node if the node belongs to the output set of the arc (i.e., if $\mathbb{T}_{va}>0$). 

The nodes participating in the self-amplifying structure are highlighted in color (orange for inputs and blue for outputs). Note that the two solutions differ in the number of nodes involved in the structure. Although every synergistic self-amplifying subhypergraph is also self-amplifying, the converse does not necessarily hold, since the Synergistic Flow Law imposes more restrictive conditions.

In the absence of interaction constraints, multiple hyperarcs can be simultaneously active, leading to a richer but less coordinated production pattern. When interaction constraints are enforced, only a subset of mutually compatible hyperarcs is selected, resulting in a more structured and coordinated configuration. The highlighted elements in the figure indicate the active subhypergraph.

In fact, the solution obtained in the absence of the \emph{Synergistic Flow Law} (Figure \ref{f:toy_example_t0_MAK0}) becomes infeasible when this law is enforced, as there exists no flow vector $\mathbf{f}$ that simultaneously satisfies both net positive realizability and the synergistic flow law. Given the matrices $\mathbb{S}$ and $\mathbb{T}$ defined above, we have $\mathbb{Q} = \mathbb{T} - \mathbb{S}$, which is given as follows, and we also show the product between $\kappa$ and $\mathbf{f}$ for the activated reactions (R2, R3, R4, R5, R6):
$$
\mathbb{Q} =
\begin{pmatrix}
-1 & 0 & -1 & 1 & 0 & 0 \\
1 & -1 & 0 & 0 & 1 & 1 \\
0 & 0 & 1 & -1 & -1 & -1 \\
1 & 0 & 0 & 0 & 0 & 0 \\
0 & 0 & 0 & 1 & 0 & -1 \\
0 & 0 & 0 & 0 & 0 & 1 \\
0 & 0 & 0 & 0 & -1 & 0
\end{pmatrix},
\qquad
\mathbf{\kappa f} = 
\begin{pmatrix}
0 \\
4.4 x_6 \\
6.7 x_3 x_6 \\
4 x_5\\
1.5 x_2 \\
6.8 x_2 x_6
\end{pmatrix}.
$$

The conditions ensuring positive flows (see equation~\eqref{ctr:12_0}) are for each $v \in \mathcal{M} = \{\text{S1, S2, S3, S5, S6}\}$:
\begin{align*}
    -6.7x_3x_6 + 4 x_5 &> 0,\\
    4.4x_6 -6.7x_3x_6 + 1.5x_2 + 6.8x_2x_6 & > 0,\\
    6.7x_3x_6 - 4x_5 -1.5x_2 -6.8x_2x_6 & > 0,\\
    4x_5 -6.8x_2x_6 & > 0,\\
    6.8x_2x_6 & > 0.
\end{align*}

Note that the system above is infeasible. Adding up the first and third inequalities we got that $-1.5x_2 - 6.8x_2x_6 > 0$. Since $x_2>0$ and $x_6>0$, it follows that $1.5 + 6.8x_6 > 0$, and then $-x_2(1.5 + 6.8x_6) < 0$, contradicting the inequality. Thus, in general,  a self-amplifying subhypergraph is not synergistic self-amplifying subhypergraph for the same input hypergraph.

\begin{figure}[htbp]
\centering
\fbox{
\begin{subfigure}{0.32\textwidth}

\resizebox{\linewidth}{!}{\begin{tikzpicture}[>=Latex, line cap=round, line join=round, x=2cm, y=1.2cm]

% ---- estilos ----
\tikzset{
  species/.style={draw, circle, minimum size=6mm, inner sep=1pt, font=\scriptsize},
  species-inactive/.style={species, fill=white},
  species-oldactive/.style={species, fill=orange!20},
  species-newactive/.style={species, fill=orange!70},
  species-out-oldactive/.style={species, fill=blue!20},
  species-out-newactive/.style={species, fill=blue!70},
  reaction-active/.style={draw=red, rectangle, minimum size=5mm, inner sep=2pt, font=\scriptsize, line width=0.8pt},
  reaction-inactive/.style={draw=black, rectangle, minimum size=5mm, inner sep=2pt, font=\scriptsize, fill=white},
  species-hidden/.style={species, draw=none, fill=none},
  reaction-hidden/.style={rectangle, minimum size=5mm, inner sep=2pt, draw=none, fill=none},
  edge-in-active/.style={->, draw=orange!80, line width=0.8pt},
  edge-in-inactive/.style={->, draw=orange!30, line width=0.6pt},
  edge-out-active/.style={->, draw=blue!90, line width=0.8pt},
  edge-out-inactive/.style={->, draw=blue!30, line width=0.6pt}
}

% --------------------
% NODOS IZQUIERDA
% --------------------
\node[species-newactive] (S0L) at (0,6) {S1};
\node[species-hidden] (S1L) at (0,5) {};
\node[species-inactive] (S2L) at (0,4) {S3};
\node[species-newactive] (S3L) at (0,3) {S4};
\node[species-hidden] (S4L) at (0,2) {};
\node[species-hidden] (S5L) at (0,1) {};
\node[species-hidden] (S6L) at (0,0) {};

% --------------------
% REACCIONES
% --------------------
\node[reaction-active] (R0) at (2,5.5) {R1};
\node[reaction-hidden] (R1) at (2,4.5) {};
\node[reaction-hidden] (R2) at (2,3.5) {};
\node[reaction-active] (R3) at (2,2.5) {R4};
\node[reaction-active] (R4) at (2,1.5) {R5};
\node[reaction-hidden] (R5) at (2,0.5) {};

% --------------------
% NODOS DERECHA
% --------------------
\node[species-out-newactive] (S0R) at (4,6) {S1};
\node[species-inactive] (S1R) at (4,5) {S2};
\node[species-hidden] (S2R) at (4,4) {};
\node[species-out-newactive] (S3R) at (4,3) {S4};
\node[species-hidden] (S4R) at (4,2) {};
\node[species-hidden] (S5R) at (4,1) {};
\node[species-inactive] (S6R) at (4,0) {S7};

% =========================================================
% ARISTAS izquierda -> reacción
% =========================================================
\draw[edge-in-active] (S3L) -- (R0);
\draw[edge-in-active] (S2L) -- (R3);
\draw[edge-in-active] (S0L) -- (R3);
\draw[edge-in-active] (S0L) -- (R4);
\draw[edge-in-active] (S3L) -- (R4);

% =========================================================
% ARISTAS reacción -> derecha
% =========================================================
\draw[edge-out-active] (R0) -- (S1R);
\draw[edge-out-active] (R0) -- (S0R);

\draw[edge-out-active] (R3) -- (S0R);

\draw[edge-out-active] (R4) -- (S3R);
\draw[edge-out-active] (R4) -- (S6R);

\end{tikzpicture}}
\caption{$t=1$.}
\end{subfigure}
\begin{subfigure}{0.32\textwidth}
\resizebox{\linewidth}{!}{\begin{tikzpicture}[>=Latex, line cap=round, line join=round, x=2cm, y=1.2cm]

% ---- estilos ----
\tikzset{
  species/.style={draw, circle, minimum size=6mm, inner sep=1pt, font=\scriptsize},
  species-inactive/.style={species, fill=white},
  species-oldactive/.style={species, fill=orange!20},
  species-newactive/.style={species, fill=orange!70},
  species-out-oldactive/.style={species, fill=blue!20},
  species-out-newactive/.style={species, fill=blue!70},
  reaction-active/.style={draw=red, rectangle, minimum size=5mm, inner sep=2pt, font=\scriptsize, line width=0.8pt},
  reaction-inactive/.style={draw=black, rectangle, minimum size=5mm, inner sep=2pt, font=\scriptsize, fill=white},
  species-hidden/.style={species, draw=none, fill=none},
  reaction-hidden/.style={rectangle, minimum size=5mm, inner sep=2pt, draw=none, fill=none},
  edge-in-active/.style={->, draw=orange!80, line width=0.8pt},
  edge-in-inactive/.style={->, draw=orange!30, line width=0.6pt},
  edge-out-active/.style={->, draw=blue!90, line width=0.8pt},
  edge-out-inactive/.style={->, draw=blue!30, line width=0.6pt}
}

% --------------------
% NODOS IZQUIERDA
% --------------------
\node[species-oldactive] (S0L) at (0,6) {S1};
\node[species-hidden] (S1L) at (0,5) {};
\node[species-inactive] (S2L) at (0,4) {S3};
\node[species-oldactive] (S3L) at (0,3) {S4};
\node[species-hidden] (S4L) at (0,2) {};
\node[species-inactive] (S5L) at (0,1) {S6};
\node[species-hidden] (S6L) at (0,0) {};
% --------------------
% REACCIONES
% --------------------
\node[reaction-inactive] (R0) at (2,5.5) {R1};
\node[reaction-hidden] (R1) at (2,4.5) {};
\node[reaction-hidden] (R2) at (2,3.5) {};
\node[reaction-inactive] (R3) at (2,2.5) {R4};
\node[reaction-inactive] (R4) at (2,1.5) {R5};
\node[reaction-active] (R5) at (2,0.5) {R6};

% --------------------
% NODOS DERECHA
% --------------------
\node[species-out-oldactive] (S0R) at (4,6) {S1};
\node[species-inactive] (S1R) at (4,5) {S2};
\node[species-hidden] (S2R) at (4,4) {};
\node[species-out-oldactive] (S3R) at (4,3) {S4};
\node[species-hidden] (S4R) at (4,2) {};
\node[species-hidden] (S5R) at (4,1) {};
\node[species-inactive] (S6R) at (4,0) {S7};

% =========================================================
% ARISTAS izquierda -> reacción
% =========================================================
\draw[edge-in-inactive] (S3L) -- (R0);
\draw[edge-in-inactive] (S2L) -- (R3);
\draw[edge-in-inactive] (S0L) -- (R3);
\draw[edge-in-inactive] (S0L) -- (R4);
\draw[edge-in-inactive] (S3L) -- (R4);

\draw[edge-in-active] (S0L) -- (R5);
\draw[edge-in-active] (S5L) -- (R5);

% =========================================================
% ARISTAS reacción -> derecha
% =========================================================
\draw[edge-out-inactive] (R0) -- (S1R);
\draw[edge-out-inactive] (R0) -- (S0R);

\draw[edge-out-inactive] (R3) -- (S0R);

\draw[edge-out-inactive] (R4) -- (S3R);
\draw[edge-out-inactive] (R4) -- (S6R);

\draw[edge-out-active] (R5) -- (S3R);
\draw[edge-out-active] (R5) -- (S6R);

\end{tikzpicture}}
\caption{$t=2$.}
\end{subfigure}
\begin{subfigure}{0.32\textwidth}
\resizebox{\linewidth}{!}{\begin{tikzpicture}[>=Latex, line cap=round, line join=round, x=2cm, y=1.2cm]

% ---- estilos ----
\tikzset{
  species/.style={draw, circle, minimum size=6mm, inner sep=1pt, font=\scriptsize},
  species-inactive/.style={species, fill=white},
  species-oldactive/.style={species, fill=orange!20},
  species-newactive/.style={species, fill=orange!70},
  species-out-oldactive/.style={species, fill=blue!20},
  species-out-newactive/.style={species, fill=blue!70},
  reaction-active/.style={draw=red, rectangle, minimum size=5mm, inner sep=2pt, font=\scriptsize, line width=0.8pt},
  reaction-inactive/.style={draw=black, rectangle, minimum size=5mm, inner sep=2pt, font=\scriptsize, fill=white},
  species-hidden/.style={species, draw=none, fill=none},
  reaction-hidden/.style={rectangle, minimum size=5mm, inner sep=2pt, draw=none, fill=none},
  edge-in-active/.style={->, draw=orange!80, line width=0.8pt},
  edge-in-inactive/.style={->, draw=orange!30, line width=0.6pt},
  edge-out-active/.style={->, draw=blue!90, line width=0.8pt},
  edge-out-inactive/.style={->, draw=blue!30, line width=0.6pt}
}

% --------------------
% NODOS IZQUIERDA
% --------------------
\node[species-oldactive] (S0L) at (0,6) {S1};
\node[species-newactive] (S1L) at (0,5) {S2};
\node[species-inactive] (S2L) at (0,4) {S3};
\node[species-oldactive] (S3L) at (0,3) {S4};
\node[species-hidden] (S4L) at (0,2) {};
\node[species-inactive] (S5L) at (0,1) {S6};
\node[species-hidden] (S6L) at (0,0) {};

% --------------------
% REACCIONES
% --------------------
\node[reaction-inactive] (R0) at (2,5.5) {R1};
\node[reaction-hidden] (R1) at (2,4.5) {};
\node[reaction-active] (R2) at (2,3.5) {R3};
\node[reaction-inactive] (R3) at (2,2.5) {R4};
\node[reaction-inactive] (R4) at (2,1.5) {R5};
\node[reaction-inactive] (R5) at (2,0.5) {R6};
% --------------------
% NODOS DERECHA
% --------------------
\node[species-out-oldactive] (S0R) at (4,6) {S1};
\node[species-out-newactive] (S1R) at (4,5) {S2};
\node[species-hidden] (S2R) at (4,4) {};
\node[species-out-oldactive] (S3R) at (4,3) {S4};
\node[species-hidden] (S4R) at (4,2) {};
\node[species-hidden] (S5R) at (4,1) {};
\node[species-inactive] (S6R) at (4,0) {S7};

% =========================================================
% ARISTAS izquierda -> reacción
% activas: S1->R2, S5->R2
% inactivas: resto
% =========================================================
\draw[edge-in-inactive] (S3L) -- (R0);

\draw[edge-in-active]   (S1L) -- (R2);
\draw[edge-in-active]   (S5L) -- (R2);

\draw[edge-in-inactive] (S2L) -- (R3);
\draw[edge-in-inactive] (S0L) -- (R3);

\draw[edge-in-inactive] (S3L) -- (R4);
\draw[edge-in-inactive] (S0L) -- (R4);

\draw[edge-in-inactive] (S5L) -- (R5);
\draw[edge-in-inactive] (S0L) -- (R5);

% =========================================================
% ARISTAS reacción -> derecha
% activas: R2->S0, R2->S6
% inactivas: resto
% =========================================================
\draw[edge-out-inactive] (R0) -- (S1R);
\draw[edge-out-inactive] (R0) -- (S0R);

\draw[edge-out-active]   (R2) -- (S0R);
\draw[edge-out-active]   (R2) -- (S6R);

\draw[edge-out-inactive] (R3) -- (S0R);

\draw[edge-out-inactive] (R4) -- (S6R);
\draw[edge-out-inactive] (R4) -- (S3R);

\draw[edge-out-inactive] (R5) -- (S6R);
\draw[edge-out-inactive] (R5) -- (S3R);

\end{tikzpicture}}
\caption{$t=3$.}
\end{subfigure}
}
\caption{A nested self-amplifying chain of subhypergraphs for $3$ periods.}
\label{f:toy_example_t3}
\end{figure}
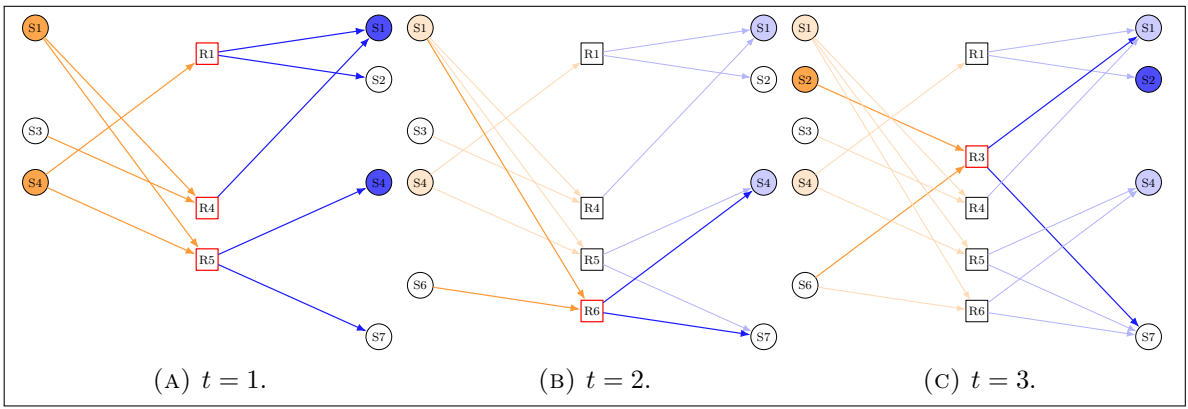

In Figure \ref{f:toy_example_t3}, we illustrate the structure of a solution for the multiperiod setting (with $t=3$). In this case, the activation decisions and state variables evolve over time. In the first period, the solution activates the self-amplifying nodes $S1$ and $S4$, together with hyperarcs $R1$, $R4$, and $R5$, using additional nodes $S3$ (input), and $S2$ and $S7$ (output). In the second period, no new self-amplifying nodes are activated, but hyperarc $R6$ is introduced, along with the extra input node $S6$, which is required to activate $R6$. In the final period, node $S2$ becomes self-amplifying, and hyperarc $R3$ is also activated.

The multiperiod dynamics reveal how the system progressively reorganizes its active structure. Early periods are driven by initially available nodes, which enable a limited set of hyperarcs. As outputs accumulate, new nodes become active, allowing additional hyperarcs to be triggered in subsequent periods. This generates a path dependent expansion process.

In particular, the figures show how certain hyperarcs emerge only after the activation of prerequisite nodes, illustrating the role of complementarities and synergies. Over time, the system converges toward a self-sustaining configuration in which activated nodes support further production, capturing the formation of a self-amplifying structure.

This example highlights how GEM captures both combinatorial selection effects and dynamic amplification, providing a unified view of structural growth in networked systems.
\section{Multiperiod Growth under GEM with Ordered Objectives}\label{sec:2}

In this section, we study the construction of expanding self-amplifying hypergraphs over multiple periods within the GEM framework through an equitable formulation, denoted as \emph{GEM-E}. The focus is on how a system evolves over time when structural feasibility, complementarities, and self-reinforcing mechanisms are jointly enforced. From an economic perspective, this corresponds to designing growth paths in which new activities emerge progressively while maintaining internal consistency and sustainability.

Rather than emphasizing extreme outcomes across periods, GEM-E adopts an ordered objective that captures the distribution of growth over time. This allows us to model preferences for balanced yet flexible expansion patterns, reflecting trade-offs between robustness, equity, and performance that are central in economic and operational settings. Beyond its intrinsic interest, this formulation provides the structural basis for the dynamic extension developed in the next section, where we introduce \emph{GEM-D}, a flow-based model that incorporates interaction-driven dynamics to capture how growth propagates endogenously across periods.

To formalize the expansion process, we introduce binary variables that capture the activation of hyperarcs over time. For each hyperarc $a \in \A$ and period $t \in \mathcal{T}$, let
$$
z_a^t =
\begin{cases}
1, & \text{if hyperarc } a \text{ is active at period } t,\\
0, & \text{otherwise.}
\end{cases}
$$
These variables encode the evolving structure of the system, as they determine which interactions or production activities are operational at each stage.

The total number of active hyperarcs at period $t$ is given by $\sum_{a \in \A} z_a^t$. Since the model enforces monotonicity of activation across periods, this quantity is non-decreasing in $t$ and therefore captures the cumulative expansion of the system over time.

This cumulative representation naturally induces a notion of incremental growth, obtained by measuring the increase in active hyperarcs between consecutive periods. Accordingly, we quantify the growth achieved at each period through the number of newly activated hyperarcs.

To this end, we consider a multiperiod GEM formulation in which period-wise growth is explicitly defined as
$$
\theta_t =
\begin{cases}
\sum_{a \in \A} z_a^1, & t = 1,\\
\sum_{a \in \A} z_a^t - \sum_{a \in \A} z_a^{t-1}, & t = 2,\ldots,T,
\end{cases}
$$
which represents the number of hyperarcs activated for the first time at period $t$. These quantities capture the marginal contribution of each period to the overall expansion and provide a direct measure of the system’s growth dynamics.

To evaluate the temporal profile of growth, we use an ordered aggregation of the vector $(\theta_1,\ldots,\theta_T)$. Instead of emphasizing peak expansion phases, we focus on the weakest periods of growth, as these are critical for ensuring a balanced and sustainable development path. In particular, we consider the sum of the $q$ smallest values of $(\theta_1,\ldots,\theta_T)$, which promotes uniform expansion across periods by penalizing low-activity periods.

\begin{definition}[Ordered $q$-sum objective]
Let $\boldsymbol{\theta} = (\theta_1,\ldots,\theta_T) \in \mathbb{R}^T$. Let $\theta_{(1)} \leq \theta_{(2)} \leq \cdots \leq \theta_{(T)}$ denote the components of $\boldsymbol{\theta}$ sorted in non-decreasing order. For a given parameter $q \in \{1,\ldots,T\}$, the ordered $q$-sum is defined as
$$
\sum_{t=1}^q \theta_{(t)}.
$$
\end{definition}

This objective belongs to the family of ordered weighted aggregation operators and provides a flexible mechanism to model equity and robustness in growth processes. In particular, $q=1$ corresponds to maximizing the minimum growth level across periods, recovering a max–min criterion. At the other extreme, $q=T$ captures total cumulative growth. Intermediate values of $q$ balance these perspectives by ensuring that a subset of the weakest periods achieves sufficiently strong expansion, without requiring full uniformity.

Given a directed multi-hypergraph $\mathcal{H}=(\mathcal{N}, \A)$ and a time horizon $T$, the multiperiod GEM-E problem can then be stated as:
\begin{align}
\max_{(\mathcal{H}^1;\mathcal{M}^1), \ldots, (\mathcal{H}^T;\mathcal{M}^T)} 
& \quad \sum_{t=1}^q \theta_{(t)} \tag{GEM-E}\label{geme}\\
\mbox{s.t.} 
& \quad \{(\mathcal{H}^t;\mathcal{M}^t)\}_{t\in \mathcal{T}} \text{ is a self-amplifying GEM.}\nonumber
\end{align}

This formulation captures the selection of a sequence of nested and self-amplifying subhypergraphs that jointly determine the system's growth trajectory. By maximizing the contribution of the weakest growth periods, the ordered objective explicitly promotes balanced expansion and  avoids solutions characterized by early spikes followed by periods of low activity.

To construct a tractable optimization model, we complement the hyperarc activation variables with node selection and flow variables that jointly determine feasible growth trajectories.

\begin{itemize}
    \item \textbf{Node selection variables:}
    \begin{itemize}
        \item[] $y^t_v \in \{0,1\}$: Takes value $1$ if node $v \in \mathcal{N}$ is active at period $t$, that is, if $v \in \mathcal{M}^t$, and $0$ otherwise.
    \end{itemize}
    \item \textbf{Flow variables:}
    \begin{itemize}
        \item[] $f^t_a \in \mathbb{R}_{\geq 0}$: Flow through hyperarc $a \in \A$ in period $t \in \mathcal{T}$.
    \end{itemize}
\end{itemize}

The quantities $\theta_t$ defined above are implicit functions of the $z$ variables and capture the incremental activation of hyperarcs. Since hyperarc activation is irreversible over time, these quantities correctly measure the introduction of new activities at each period and thus represent the system’s growth profile.

We now describe the constraints that characterize feasible GEM trajectories and ensure structural consistency, self-sufficiency, and effective activity of the system over time.

\begin{itemize}

\item \textbf{Monotonicity of nodes across periods:}
\begin{align}
y^t_v \leq y^{t+1}_v, \quad \forall v \in \mathcal{N}, \; t \in \mathcal{T} \backslash \{T\}. \label{ctr:3}
\end{align}
These constraints enforce persistence of activated nodes: once a node becomes active, it remains available in subsequent periods. From an economic viewpoint, this reflects the irreversibility of investments or capabilities.

\item \textbf{Monotonicity of hyperarcs across periods:}
\begin{align}
z^t_a \leq z^{t+1}_a, \quad \forall a \in \A, \; t \in \mathcal{T} \backslash \{T\}. \label{ctr:4}
\end{align}
Similarly, once a hyperarc is activated, it remains active in later periods. This captures persistence of production technologies or interactions, and guarantees that growth is cumulative.

\item \textbf{Self-sufficiency in arcs (node feasibility):}
\begin{align}
y_v^t &\leq \sum_{a\in \A \, : \, \mathbb{S}_{va} > 0} z_a^t, \quad \forall v \in \mathcal{N}, \; t\in \mathcal{T},\label{ctr:5a}\\
y_v^t &\leq \sum_{a\in \A \, : \, \mathbb{T}_{va} > 0} z_a^t, \quad \forall v \in \mathcal{N}, \; t\in \mathcal{T}.\label{ctr:5b}
\end{align}
These constraints ensure that any active node participates in at least one active hyperarc both as input and as output. In economic terms, this prevents isolated activities and enforces that each active entity is embedded within the production system.

\item \textbf{Self-sufficiency in nodes (arc feasibility):}
\begin{align}
z_a^t &\leq \sum_{v\in \mathcal{N} \, : \, \mathbb{S}_{va} > 0} y_v^t, \quad \forall a \in \A, \; t\in \mathcal{T},\label{ctr:6a}\\
z_a^t &\leq \sum_{v\in \mathcal{N} \, : \, \mathbb{T}_{va} > 0} y_v^t, \quad \forall a \in \A, \; t\in \mathcal{T}.\label{ctr:6b}
\end{align}
These constraints impose that a hyperarc can only be activated if all the nodes required for its inputs and outputs are active. This captures complementarities among activities: production processes require the simultaneous availability of all necessary components.

\item \textbf{Net positive realizability:}
\begin{align*}
y_v^t = 1 \;\Rightarrow\; \sum_{a \in \A} \mathbb{Q}_{va} f_a^t \geq \varepsilon, 
\quad \forall v \in \mathcal{N}, \; t\in \mathcal{T}.
\end{align*}
This condition ensures that each active node achieves a strictly positive net balance of flows. The above implication can be enforced through a standard big-$M$ linearization:
\begin{align}
\sum_{a \in \A} \mathbb{Q}_{va} f_a^t \geq \varepsilon - \Delta(1-y_v^t), 
\quad \forall v \in \mathcal{N}, \; t\in \mathcal{T}, \label{ctr:7}
\end{align}
\noindent where $\Delta>0$ is a sufficiently large constant. Note that when $y_v^t=1$, the constraint reduces to 
$\sum_{a \in \A} \mathbb{Q}_{va} f_a^t \geq \varepsilon$, enforcing that node $v$ is sustained by the system. 
When $y_v^t=0$, the constraint is relaxed via the big-$M$ term and does not impose any restriction. 

From an economic perspective, this condition guarantees viability: any activated node must generate a strictly positive net contribution, preventing the selection of inactive or non-productive components.

\item \textbf{Flow activation constraints:}
\begin{align}
\epsilon_a z_a \leq f_a^t &\leq \Delta_a z_a^t, \quad \forall a \in \A, \; t\in \mathcal{T},\label{ctr:8a}
\end{align}
These constraints link structural activation and operational activity. If a hyperarc is inactive, no flow can pass through it. If it is active, it must carry a strictly positive flow. Thus, only economically meaningful or operationally effective interactions are selected.
\end{itemize}

With these components, the multiperiod \eqref{geme} model consists of maximizing the ordered $q$-sum of the smallest components of the growth vector $(\theta_1,\ldots,\theta_T)$ subject to constraints \eqref{ctr:3}--\eqref{ctr:8a} and the domain conditions:
\begin{align}
& y^t_{v}\in \{0,1\},\quad \forall v \in \N, t \in \mathcal{T}, \label{ctr:MESHA_y_domain}\\
& z^t_{a}\in \{0,1\},\quad \forall a \in \A, t \in \mathcal{T},\label{ctr:MESHA_z_domain}\\
& f^t_a \in \mathbb{R}_{>0},\quad \forall  a \in \A, t \in \mathcal{T}.\label{ctr:MESHA_f_domain}
\end{align}

It is important to note that the objective function defined above is not linear, as it involves the ordered selection of the smallest components of the vector $(\theta_1,\ldots,\theta_T)$. This type of objective belongs to the class of ordered weighted aggregation functions, which are widely used to model equity and robustness considerations in optimization.

To obtain a tractable formulation, we can rely on equivalent linear representations of ordered objectives available in the literature. In particular, the problem of minimizing or maximizing the sum of the $q$ smallest (or largest) components can be reformulated using auxiliary variables and linear constraints. Classical results by \citet{ogryczak2003minimizing} provide compact formulations for such ordered sums. Alternatively, one can use the representation proposed in \citet{cont:blanco2014revisiting}, which has proven effective in continuous and discrete optimization settings. Notably, the linear relaxation of the latter coincides with that of the former, while offering a modeling structure that is particularly convenient for integration within mixed-integer programming frameworks.

In our implementation, we adopt the formulation in \citet{cont:blanco2014revisiting}, which enables us to embed the ordered $q$-sum objective within a mixed-integer linear optimization model without loss of tractability. This formulation is:
\begin{align*}
    \max\;&\sum_{t \in \mathcal{T}} u_t + \sum_{k\in \mathcal{T}} w_k\\
    \mbox{s.t. }  
    & \eqref{ctr:3}-\eqref{ctr:8a},\\
    & u_t + w_k \leq \theta_t, \qquad \forall v \in \mathcal{T}, k=1, \ldots, q,\\
    & u_v + w_k \leq 0,\qquad  \forall k=q+1, \ldots, T,\;  v \in \mathcal{T},\\
    & u_t \in \mathbb{R},\qquad  \forall t \in \mathcal{T},\\
    & w_k \in \mathbb{R},\qquad  \forall k \in \mathcal{T},\\
&\eqref{ctr:MESHA_y_domain},\eqref{ctr:MESHA_z_domain}, \eqref{ctr:MESHA_f_domain}.
\end{align*}

From an economic viewpoint, this formulation represents a dynamic growth planning problem in which complementarities among activities shape feasible expansion paths. By prioritizing the weakest periods of growth, the ordered objective promotes balanced development, mitigates bottlenecks, and ensures that expansion is sustained over time rather than concentrated in a few periods. This aligns naturally with economic notions of coordinated development and stable growth trajectories.
\section{Flow Synergistic Multiperiod GEM Problems}\label{sec:3}

The structural formulation developed in the previous section establishes the foundations for the expansion of self-amplifying hypergraphs over multiple periods within the GEM framework. However, that model captures only the combinatorial and topological features of the evolution process, without explicitly incorporating the dynamic interactions among nodes and hyperarcs that drive the activation of new elements. 

In this section, we extend the GEM formulation by integrating flow-dependent dynamics governed by the \emph{Synergistic Flow Law}, which introduces nonlinear relations between node activities and hyperarc intensities. This extension enables a richer representation of growth processes, where activation decisions are not only structurally feasible but also dynamically consistent with the underlying interaction mechanisms.

The resulting model, referred to as the \textit{Flow Synergistic Multiperiod GEM Problem (GEM-D)}, combines the structural consistency of the GEM formulation with nonlinear flow constraints. This leads to a coupled system in which structure and dynamics co-evolve over time, capturing how complementarities among active nodes generate and sustain new interactions.

Specifically,  the GEM-D model is a restricted version of \eqref{geme} that can be stated as:
\begin{align}
\max_{(\mathcal{H}^1;\mathcal{M}^1), \ldots, (\mathcal{H}^T;\mathcal{M}^T)} 
& \quad \sum_{t=1}^q \theta_{(t)} \tag{GEM-D}\label{gemd}\\
\mbox{s.t.} 
& \quad \{(\mathcal{H}^t;\mathcal{M}^t)\}_{t\in \mathcal{T}} \text{ is a self-amplifying GEM,}\nonumber\\
& \quad \text{Flows at period $t$ satisfy Synergistic Flow Law \eqref{ctr:12_0}.} 
\end{align}

To derive a model for this problem, in addition to the $z$, $y$, and $f$ variables introduced in the previous section, we consider the following state variables:
$$
x_v^t \in \R_{\geq 0}: \text{ State of node $v\in \mathcal{N}$ in period $t\in \mathcal{T}$},
$$
which represent the \emph{activity level}, \emph{availability}, or \emph{intensity} associated with each node.

The state value of each node evolves over time as a function of its previous state and the flows generated in the current period. To initialize the system, we introduce variables $x_v^0 \in \mathbb{R}_{\geq 0}$ representing the \emph{initial states}, i.e., to represent the required initial amount of the nodes to start the growing process. 

The state evolution is given by:
\begin{align}
x_v^t = x_v^{t-1} + \sum_{a \in \A} \mathbb{Q}_{va} f_a^{t}, 
\quad \forall v \in \mathcal{N}, \; t \in \mathcal{T} (t>1), 
\label{ctr:9}
\end{align}
where the state of node $v$ at period $t$ results from its previous state plus the net contribution of flows. 

To explicitly control node activation, we introduce binary variables:
\begin{align}
\rho_v^t \in \{0,1\}: \text{ 1 if node $v$ is active in period $t$, and 0 otherwise}, \quad \forall t \in \mathcal{T}.
\end{align}

The nested structure is enforced through:
\begin{align}
\rho_v^t \leq \rho_v^{t+1}, 
\quad \forall v \in \mathcal{N}, \; t \in \mathcal{T} \setminus \{T\},
\label{ctr:10}
\end{align}
ensuring persistence of active nodes.

State variables are linked to activation through:
\begin{align}
\epsilon \, \rho_v^t \leq x_v^t \leq \frac{1}{\epsilon} \rho_v^t, 
\quad \forall v \in \mathcal{N}, \; t \in \mathcal{T}, 
\label{ctr:11}
\end{align}
so that inactive nodes have zero state, while active nodes have strictly positive values.

The \emph{Synergistic Flow Law} is modeled as:
\begin{equation}\label{ctr:12}
   f_a^t = \kappa_a z_a^t \prod_{v \in \mathcal{N}} (x_v^t)^{\mathbb{S}_{va}}, 
\forall a \in \A, \; t \in \mathcal{T}
\end{equation}

In our setting, the state variables $x_v^t$ play an analogous role to concentrations, while the stoichiometric coefficients $\mathbb{S}_{va}$ and $\mathbb{T}_{va}$ determine the influence of each node on the activation of hyperarc $a$.

A key distinction with standard mass-action models is the presence of the activation variable $z_a^t$, which endogenously determines whether the hyperarc is active in the network. This introduces a coupling between discrete design decisions and continuous flow dynamics: even if the state variables would allow a positive flow, the hyperarc must be selected in the evolving subhypergraph for the interaction to take place.

From an economic and operational perspective, this law captures how cooperative interactions among active components generate productive flows. The multiplicative structure encodes synergy effects, whereby the joint presence of several nodes amplifies activity, while the reversibility mechanism allows for feedback and dissipation within the system.

Combining all these elements, \eqref{gemd} can be formulated as the following Mixed-Integer Nonlinear Programming (MINLP) model:
\begin{align}
\max_{x,\,f,\,y,\,z,\,\rho} \; & \sum_{t=1}^q \theta_{(t)} \nonumber\\[3pt]
\text{s.t.} \quad
    & \eqref{ctr:3}- \eqref{ctr:8a}, \eqref{ctr:9}-  \eqref{ctr:12}, \nonumber\\
    & y^t_{v} \in \{0,1\}, \quad \forall v \in \mathcal{N}, \; t \in \mathcal{T}, \\
    & z^t_{a} \in \{0,1\}, \quad \forall a \in \A, \; t \in \mathcal{T}, \\
    & f^t_a \in \mathbb{R}_{+}, \quad \forall a \in \A, \; t \in \mathcal{T}, \\
    & x^t_v \in \mathbb{R}_{+}, \quad \forall v \in \mathcal{N}, \; t \in \mathcal{T}, \\
    & \rho^t_v \in \{0,1\}, \quad \forall v \in \mathcal{N}, \; t \in \mathcal{T}.
\end{align}

This formulation extends the structural GEM model by enforcing that growth is not only combinatorially feasible but also dynamically consistent with synergistic interactions. As a result, the expansion process must satisfy both structural and flow-based constraints. The reformulation of the highly nonlinear constraints \eqref{ctr:12} as logarithmic constraints that can be handled by the solvers is detailed in \ref{app:1}

From a computational perspective, this approach leads to a mixed-integer linear programming (MILP) formulation that can be efficiently handled by modern solvers such as \texttt{Gurobi} or \texttt{CPLEX}. Moreover, the structure of the approximation allows for a favorable trade-off between numerical stability and solution quality, making it suitable for large-scale instances of \eqref{gemd}.

\begin{figure}[htbp]
\centering
\fbox{
\begin{subfigure}{0.32\textwidth}

\resizebox{\linewidth}{!}{\begin{tikzpicture}[>=Latex, line cap=round, line join=round, x=2cm, y=1.2cm]

% ---- estilos ----
\tikzset{
  species/.style={draw, circle, minimum size=6mm, inner sep=1pt, font=\scriptsize},
  species-inactive/.style={species, fill=white},
  species-oldactive/.style={species, fill=orange!20},
  species-newactive/.style={species, fill=orange!70},
  species-out-oldactive/.style={species, fill=blue!20},
  species-out-newactive/.style={species, fill=blue!70},
  reaction-active/.style={draw=red, rectangle, minimum size=5mm, inner sep=2pt, font=\scriptsize, line width=0.8pt},
  reaction-inactive/.style={draw=black, rectangle, minimum size=5mm, inner sep=2pt, font=\scriptsize, fill=white},
  species-hidden/.style={species, draw=none, fill=none},
  reaction-hidden/.style={rectangle, minimum size=5mm, inner sep=2pt, draw=none, fill=none},
  edge-in-active/.style={->, draw=orange!80, line width=0.8pt},
  edge-in-inactive/.style={->, draw=orange!30, line width=0.6pt},
  edge-out-active/.style={->, draw=blue!90, line width=0.8pt},
  edge-out-inactive/.style={->, draw=blue!30, line width=0.6pt}
}

% --------------------
% NODOS IZQUIERDA
% --------------------
\node[species-newactive] (S0L) at (0,6) {S1};
\node[species-hidden] (S1L) at (0,5) {};
\node[species-inactive] (S2L) at (0,4) {S3};
\node[species-newactive] (S3L) at (0,3) {S4};
\node[species-hidden] (S4L) at (0,2) {};
\node[species-hidden] (S5L) at (0,1) {};
\node[species-hidden] (S6L) at (0,0) {};

% --------------------
% REACCIONES
% --------------------
\node[reaction-active] (R0) at (2,5.5) {R1};
\node[reaction-hidden] (R1) at (2,4.5) {};
\node[reaction-hidden] (R2) at (2,3.5) {};
\node[reaction-active] (R3) at (2,2.5) {R4};
\node[reaction-active] (R4) at (2,1.5) {R5};
\node[reaction-hidden] (R5) at (2,0.5) {};

% --------------------
% NODOS DERECHA
% --------------------
\node[species-out-newactive] (S0R) at (4,6) {S1};
\node[species-inactive] (S1R) at (4,5) {S2};
\node[species-hidden] (S2R) at (4,4) {};
\node[species-out-newactive] (S3R) at (4,3) {S4};
\node[species-hidden] (S4R) at (4,2) {};
\node[species-hidden] (S5R) at (4,1) {};
\node[species-inactive] (S6R) at (4,0) {S7};

% =========================================================
% ARISTAS izquierda -> reacción
% =========================================================
\draw[edge-in-active] (S3L) -- (R0);
\draw[edge-in-active] (S2L) -- (R3);
\draw[edge-in-active] (S0L) -- (R3);
\draw[edge-in-active] (S0L) -- (R4);
\draw[edge-in-active] (S3L) -- (R4);

% =========================================================
% ARISTAS reacción -> derecha
% =========================================================
\draw[edge-out-active] (R0) -- (S1R);
\draw[edge-out-active] (R0) -- (S0R);

\draw[edge-out-active] (R3) -- (S0R);

\draw[edge-out-active] (R4) -- (S3R);
\draw[edge-out-active] (R4) -- (S6R);

\end{tikzpicture}}
\caption{$t=1$}
\end{subfigure}
\begin{subfigure}{0.32\textwidth}
\resizebox{\linewidth}{!}{\begin{tikzpicture}[>=Latex, line cap=round, line join=round, x=2cm, y=1.2cm]

% ---- estilos ----
\tikzset{
  species/.style={draw, circle, minimum size=6mm, inner sep=1pt, font=\scriptsize},
  species-inactive/.style={species, fill=white},
  species-oldactive/.style={species, fill=orange!20},
  species-newactive/.style={species, fill=orange!70},
  species-out-oldactive/.style={species, fill=blue!20},
  species-out-newactive/.style={species, fill=blue!70},
  reaction-active/.style={draw=red, rectangle, minimum size=5mm, inner sep=2pt, font=\scriptsize, line width=0.8pt},
  reaction-inactive/.style={draw=black, rectangle, minimum size=5mm, inner sep=2pt, font=\scriptsize, fill=white},
  species-hidden/.style={species, draw=none, fill=none},
  reaction-hidden/.style={rectangle, minimum size=5mm, inner sep=2pt, draw=none, fill=none},
  edge-in-active/.style={->, draw=orange!80, line width=0.8pt},
  edge-in-inactive/.style={->, draw=orange!30, line width=0.6pt},
  edge-out-active/.style={->, draw=blue!90, line width=0.8pt},
  edge-out-inactive/.style={->, draw=blue!30, line width=0.6pt}
}

% --------------------
% NODOS IZQUIERDA
% --------------------
\node[species-oldactive] (S0L) at (0,6) {S1};
\node[species-hidden] (S1L) at (0,5) {};
\node[species-inactive] (S2L) at (0,4) {S3};
\node[species-oldactive] (S3L) at (0,3) {S4};
\node[species-hidden] (S4L) at (0,2) {};
\node[species-inactive] (S5L) at (0,1) {S6};
\node[species-hidden] (S6L) at (0,0) {};
% --------------------
% REACCIONES
% --------------------
\node[reaction-inactive] (R0) at (2,5.5) {R1};
\node[reaction-hidden] (R1) at (2,4.5) {};
\node[reaction-hidden] (R2) at (2,3.5) {};
\node[reaction-inactive] (R3) at (2,2.5) {R4};
\node[reaction-inactive] (R4) at (2,1.5) {R5};
\node[reaction-active] (R5) at (2,0.5) {R6};

% --------------------
% NODOS DERECHA
% --------------------
\node[species-out-oldactive] (S0R) at (4,6) {S1};
\node[species-inactive] (S1R) at (4,5) {S2};
\node[species-hidden] (S2R) at (4,4) {};
\node[species-out-oldactive] (S3R) at (4,3) {S4};
\node[species-hidden] (S4R) at (4,2) {};
\node[species-hidden] (S5R) at (4,1) {};
\node[species-inactive] (S6R) at (4,0) {S7};

% =========================================================
% ARISTAS izquierda -> reacción
% =========================================================
\draw[edge-in-inactive] (S3L) -- (R0);
\draw[edge-in-inactive] (S2L) -- (R3);
\draw[edge-in-inactive] (S0L) -- (R3);
\draw[edge-in-inactive] (S0L) -- (R4);
\draw[edge-in-inactive] (S3L) -- (R4);

\draw[edge-in-active] (S0L) -- (R5);
\draw[edge-in-active] (S5L) -- (R5);

% =========================================================
% ARISTAS reacción -> derecha
% =========================================================
\draw[edge-out-inactive] (R0) -- (S1R);
\draw[edge-out-inactive] (R0) -- (S0R);

\draw[edge-out-inactive] (R3) -- (S0R);

\draw[edge-out-inactive] (R4) -- (S3R);
\draw[edge-out-inactive] (R4) -- (S6R);

\draw[edge-out-active] (R5) -- (S3R);
\draw[edge-out-active] (R5) -- (S6R);

\end{tikzpicture}}
\caption{$t=2$}
\end{subfigure}
\begin{subfigure}{0.32\textwidth}
\resizebox{\linewidth}{!}{\begin{tikzpicture}[>=Latex, line cap=round, line join=round, x=2cm, y=1.2cm]

% ---- estilos ----
\tikzset{
  species/.style={draw, circle, minimum size=6mm, inner sep=1pt, font=\scriptsize},
  species-inactive/.style={species, fill=white},
  species-oldactive/.style={species, fill=orange!20},
  species-newactive/.style={species, fill=orange!70},
  species-out-oldactive/.style={species, fill=blue!20},
  species-out-newactive/.style={species, fill=blue!70},
  reaction-active/.style={draw=red, rectangle, minimum size=5mm, inner sep=2pt, font=\scriptsize, line width=0.8pt},
  reaction-inactive/.style={draw=black, rectangle, minimum size=5mm, inner sep=2pt, font=\scriptsize, fill=white},
  species-hidden/.style={species, draw=none, fill=none},
  reaction-hidden/.style={rectangle, minimum size=5mm, inner sep=2pt, draw=none, fill=none},
  edge-in-active/.style={->, draw=orange!80, line width=0.8pt},
  edge-in-inactive/.style={->, draw=orange!30, line width=0.6pt},
  edge-out-active/.style={->, draw=blue!90, line width=0.8pt},
  edge-out-inactive/.style={->, draw=blue!30, line width=0.6pt}
}

% --------------------
% NODOS IZQUIERDA
% --------------------
\node[species-oldactive] (S0L) at (0,6) {S1};
\node[species-newactive] (S1L) at (0,5) {S2};
\node[species-inactive] (S2L) at (0,4) {S3};
\node[species-oldactive] (S3L) at (0,3) {S4};
\node[species-hidden] (S4L) at (0,2) {};
\node[species-inactive] (S5L) at (0,1) {S6};
\node[species-hidden] (S6L) at (0,0) {};

% --------------------
% REACCIONES
% --------------------
\node[reaction-inactive] (R0) at (2,5.5) {R1};
\node[reaction-hidden] (R1) at (2,4.5) {};
\node[reaction-active] (R2) at (2,3.5) {R3};
\node[reaction-inactive] (R3) at (2,2.5) {R4};
\node[reaction-inactive] (R4) at (2,1.5) {R5};
\node[reaction-inactive] (R5) at (2,0.5) {R6};
% --------------------
% NODOS DERECHA
% --------------------
\node[species-out-oldactive] (S0R) at (4,6) {S1};
\node[species-out-newactive] (S1R) at (4,5) {S2};
\node[species-hidden] (S2R) at (4,4) {};
\node[species-out-oldactive] (S3R) at (4,3) {S4};
\node[species-hidden] (S4R) at (4,2) {};
\node[species-hidden] (S5R) at (4,1) {};
\node[species-inactive] (S6R) at (4,0) {S7};

% =========================================================
% ARISTAS izquierda -> reacción
% activas: S1->R2, S5->R2
% inactivas: resto
% =========================================================
\draw[edge-in-inactive] (S3L) -- (R0);

\draw[edge-in-active]   (S1L) -- (R2);
\draw[edge-in-active]   (S5L) -- (R2);

\draw[edge-in-inactive] (S2L) -- (R3);
\draw[edge-in-inactive] (S0L) -- (R3);

\draw[edge-in-inactive] (S3L) -- (R4);
\draw[edge-in-inactive] (S0L) -- (R4);

\draw[edge-in-inactive] (S5L) -- (R5);
\draw[edge-in-inactive] (S0L) -- (R5);

% =========================================================
% ARISTAS reacción -> derecha
% activas: R2->S0, R2->S6
% inactivas: resto
% =========================================================
\draw[edge-out-inactive] (R0) -- (S1R);
\draw[edge-out-inactive] (R0) -- (S0R);

\draw[edge-out-active]   (R2) -- (S0R);
\draw[edge-out-active]   (R2) -- (S6R);

\draw[edge-out-inactive] (R3) -- (S0R);

\draw[edge-out-inactive] (R4) -- (S6R);
\draw[edge-out-inactive] (R4) -- (S3R);

\draw[edge-out-inactive] (R5) -- (S6R);
\draw[edge-out-inactive] (R5) -- (S3R);

\end{tikzpicture}}
\caption{$t=3$}
\end{subfigure}
}
\caption{Solution of \eqref{geme} for $3$ periods.}
\label{f:nomak}
\end{figure}
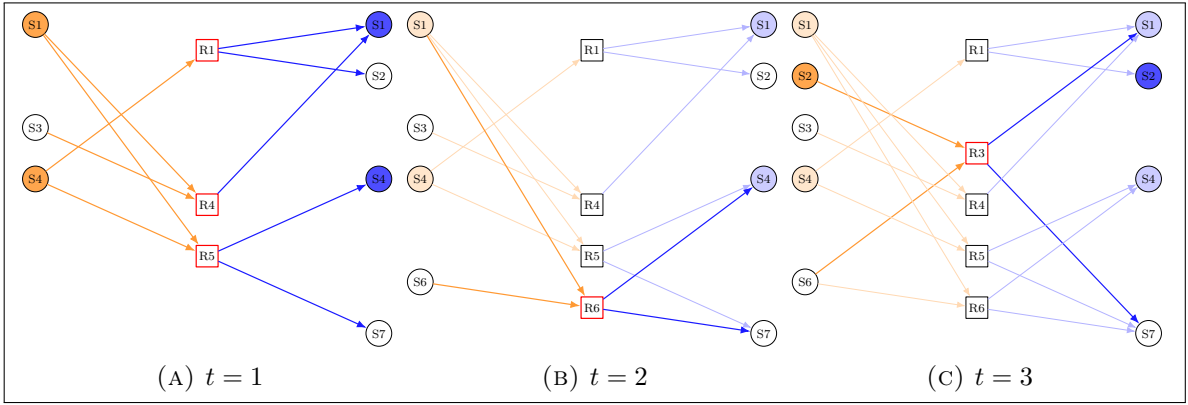

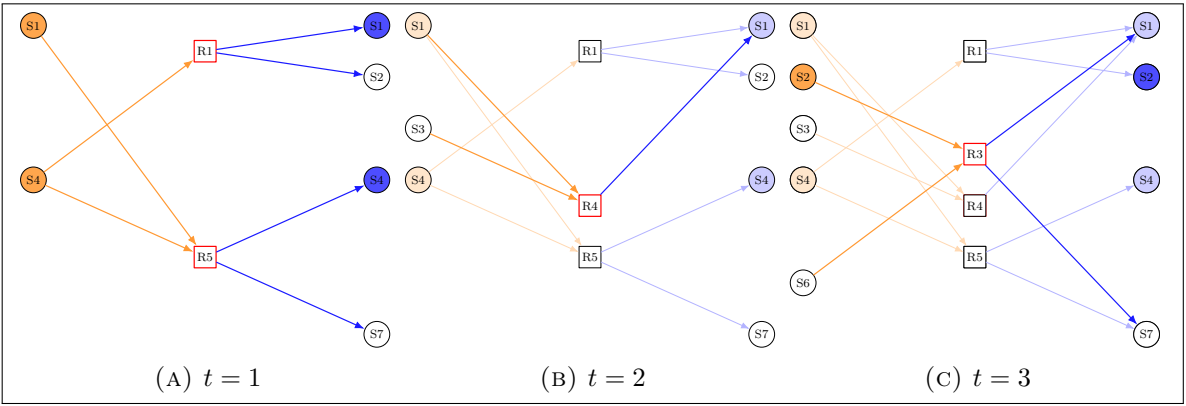
\begin{figure}[htbp]
\centering
\fbox{
\begin{subfigure}{0.32\textwidth}
\resizebox{\linewidth}{!}{\begin{tikzpicture}[>=Latex, line cap=round, line join=round, x=2cm, y=1.2cm]

\tikzset{
  species/.style={draw, circle, minimum size=6mm, inner sep=1pt, font=\scriptsize},
  species-inactive/.style={species, fill=white},
  species-oldactive/.style={species, fill=orange!20},
  species-newactive/.style={species, fill=orange!70},
  species-out-oldactive/.style={species, fill=blue!20},
  species-out-newactive/.style={species, fill=blue!70},
  reaction-active/.style={draw=red, rectangle, minimum size=5mm, inner sep=2pt, font=\scriptsize, line width=0.8pt},
  reaction-inactive/.style={draw=black, rectangle, minimum size=5mm, inner sep=2pt, font=\scriptsize, fill=white},
  species-hidden/.style={species, draw=none, fill=none},
  reaction-hidden/.style={rectangle, minimum size=5mm, inner sep=2pt, draw=none, fill=none},
  edge-in-active/.style={->, draw=orange!80, line width=0.8pt},
  edge-in-inactive/.style={->, draw=orange!30, line width=0.6pt},
  edge-out-active/.style={->, draw=blue!90, line width=0.8pt},
  edge-out-inactive/.style={->, draw=blue!30, line width=0.6pt}
}

% bounding box fija (IMPORTANTE para alinear subfiguras)
%\path[use as bounding box] (-1,-5) rectangle (11,5);

% --------------------
% NODOS IZQUIERDA
% --------------------
\node[species-newactive] (S0L) at (0,6) {S1};
\node[species-hidden] (S1L) at (0,5) {};
\node[species-hidden] (S2L) at (0,4) {};
\node[species-newactive] (S3L) at (0,3) {S4};
\node[species-hidden] (S4L) at (0,2) {};
\node[species-hidden] (S5L) at (0,1) {};
\node[species-hidden] (S6L) at (0,0) {};

% --------------------
% REACCIONES
% --------------------
\node[reaction-active] (R0) at (2,5.5) {R1};
\node[reaction-hidden] (R1) at (2,4.5) {};
\node[reaction-hidden] (R2) at (2,3.5) {};
\node[reaction-hidden] (R3) at (2,2.5) {};
\node[reaction-active] (R4) at (2,1.5) {R5};
\node[reaction-hidden] (R5) at (2,0.5) {};

% --------------------
% NODOS DERECHA
% --------------------
\node[species-out-newactive] (S0R) at (4,6) {S1};
\node[species-inactive] (S1R) at (4,5) {S2};
\node[species-hidden] (S2R) at (4,4) {};
\node[species-out-newactive] (S3R) at (4,3) {S4};
\node[species-hidden] (S4R) at (4,2) {};
\node[species-hidden] (S5R) at (4,1) {};
\node[species-inactive] (S6R) at (4,0) {S7};

% --------------------
% ARISTAS ENTRADA (azul)
% --------------------
\draw[edge-in-active] (S3L) -- (R0);
\draw[edge-in-active] (S3L) -- (R4);
\draw[edge-in-active] (S0L) -- (R4);

% --------------------
% ARISTAS SALIDA (naranja)
% --------------------
\draw[edge-out-active] (R0) -- (S0R);
\draw[edge-out-active] (R0) -- (S1R);

\draw[edge-out-active] (R4) -- (S3R);
\draw[edge-out-active] (R4) -- (S6R);

\end{tikzpicture}}
\caption{$t=1$}
\end{subfigure}
\begin{subfigure}{0.32\textwidth}
\resizebox{\linewidth}{!}{\begin{tikzpicture}[>=Latex, line cap=round, line join=round, x=2cm, y=1.2cm]

\tikzset{
  species/.style={draw, circle, minimum size=6mm, inner sep=1pt, font=\scriptsize},
  species-inactive/.style={species, fill=white},
  species-oldactive/.style={species, fill=orange!20},
  species-newactive/.style={species, fill=orange!70},
  species-out-oldactive/.style={species, fill=blue!20},
  species-out-newactive/.style={species, fill=blue!70},
  reaction-active/.style={draw=red, rectangle, minimum size=5mm, inner sep=2pt, font=\scriptsize, line width=0.8pt},
  reaction-inactive/.style={draw=black, rectangle, minimum size=5mm, inner sep=2pt, font=\scriptsize, fill=white},
  species-hidden/.style={species, draw=none, fill=none},
  reaction-hidden/.style={rectangle, minimum size=5mm, inner sep=2pt, draw=none, fill=none},
  edge-in-active/.style={->, draw=orange!80, line width=0.8pt},
  edge-in-inactive/.style={->, draw=orange!30, line width=0.6pt},
  edge-out-active/.style={->, draw=blue!90, line width=0.8pt},
  edge-out-inactive/.style={->, draw=blue!30, line width=0.6pt}
}

% bounding box fija para alinear subfiguras
%\path[use as bounding box] (-1,-5) rectangle (11,5);

% --------------------
% NODOS IZQUIERDA
% --------------------
\node[species-oldactive] (S0L) at (0,6) {S1};
\node[species-hidden] (S1L) at (0,5) {};
\node[species-inactive] (S2L) at (0,4) {S3};
\node[species-oldactive] (S3L) at (0,3) {S4};
\node[species-hidden] (S4L) at (0,2) {};
\node[species-hidden] (S5L) at (0,1) {};
\node[species-hidden] (S6L) at (0,0) {};

% --------------------
% REACCIONES
% --------------------
\node[reaction-inactive] (R0) at (2,5.5) {R1};
\node[reaction-hidden] (R1) at (2,4.5) {};
\node[reaction-hidden] (R2) at (2,3.5) {};
\node[reaction-active] (R3) at (2,2.5) {R4};
\node[reaction-inactive] (R4) at (2,1.5) {R5};
\node[reaction-hidden] (R5) at (2,0.5) {};

% --------------------
% NODOS DERECHA
% --------------------
\node[species-out-oldactive] (S0R) at (4,6) {S1};
\node[species-inactive] (S1R) at (4,5) {S2};
\node[species-hidden] (S2R) at (4,4) {};
\node[species-out-oldactive] (S3R) at (4,3) {S4};
\node[species-hidden] (S4R) at (4,2) {};
\node[species-hidden] (S5R) at (4,1) {};
\node[species-inactive] (S6R) at (4,0) {S7};

% --------------------
% ARISTAS ANTIGUAS (primero)
% --------------------
\draw[edge-in-inactive] (S3L) -- (R0);
\draw[edge-out-inactive]  (R0) -- (S1R);
\draw[edge-out-inactive]  (R0) -- (S0R);

\draw[edge-in-inactive] (S3L) -- (R4);
\draw[edge-in-inactive] (S0L) -- (R4);
\draw[edge-out-inactive]  (R4) -- (S3R);
\draw[edge-out-inactive]  (R4) -- (S6R);

% --------------------
% ARISTAS NUEVAS (después, por encima)
% --------------------
\draw[edge-in-active] (S2L) -- (R3);
\draw[edge-in-active] (S0L) -- (R3);
\draw[edge-out-active]  (R3) -- (S0R);

\end{tikzpicture}}
\caption{$t=2$}
\end{subfigure}
\begin{subfigure}{0.32\textwidth}
\resizebox{\linewidth}{!}{\begin{tikzpicture}[>=Latex, line cap=round, line join=round, x=2cm, y=1.2cm]

\tikzset{
  species/.style={draw, circle, minimum size=6mm, inner sep=1pt, font=\scriptsize},
  species-inactive/.style={species, fill=white},
  species-oldactive/.style={species, fill=orange!20},
  species-newactive/.style={species, fill=orange!70},
  species-out-oldactive/.style={species, fill=blue!20},
  species-out-newactive/.style={species, fill=blue!70},
  reaction-active/.style={draw=red, rectangle, minimum size=5mm, inner sep=2pt, font=\scriptsize, line width=0.8pt},
  reaction-inactive/.style={draw=black, rectangle, minimum size=5mm, inner sep=2pt, font=\scriptsize, fill=white},
  species-hidden/.style={species, draw=none, fill=none},
  reaction-hidden/.style={rectangle, minimum size=5mm, inner sep=2pt, draw=none, fill=none},
  edge-in-active/.style={->, draw=orange!80, line width=0.8pt},
  edge-in-inactive/.style={->, draw=orange!30, line width=0.6pt},
  edge-out-active/.style={->, draw=blue!90, line width=0.8pt},
  edge-out-inactive/.style={->, draw=blue!30, line width=0.6pt}
}

% bounding box fija para alinear subfiguras
% \path[use as bounding box] (-1,-5) rectangle (11,5);

% --------------------
% NODOS IZQUIERDA
% --------------------
\node[species-oldactive] (S0L) at (0,6) {S1};
\node[species-hidden] (S1L) at (0,5) {};
\node[species-inactive] (S2L) at (0,4) {S3};
\node[species-oldactive] (S3L) at (0,3) {S4};
\node[species-hidden] (S4L) at (0,2) {};
\node[species-hidden] (S5L) at (0,1) {};
\node[species-hidden] (S6L) at (0,0) {};

% --------------------
% REACCIONES
% --------------------
\node[reaction-inactive] (R0) at (2,5.5) {R1};
\node[reaction-hidden] (R1) at (2,4.5) {};
\node[reaction-hidden] (R2) at (2,3.5) {};
\node[reaction-active] (R3) at (2,2.5) {R4};
\node[reaction-inactive] (R4) at (2,1.5) {R5};
\node[reaction-hidden] (R5) at (2,0.5) {};

% --------------------
% NODOS DERECHA
% --------------------
\node[species-out-oldactive] (S0R) at (4,6) {S1};
\node[species-inactive] (S1R) at (4,5) {S2};
\node[species-hidden] (S2R) at (4,4) {};
\node[species-out-oldactive] (S3R) at (4,3) {S4};
\node[species-hidden] (S4R) at (4,2) {};
\node[species-hidden] (S5R) at (4,1) {};
\node[species-inactive] (S6R) at (4,0) {S7};

% --------------------
% NODOS IZQUIERDA
% --------------------
\node[species-oldactive] (S0L) at (0,6) {S1};
\node[species-newactive] (S1L) at (0,5) {S2};
\node[species-inactive] (S2L) at (0,4) {S3};
\node[species-oldactive] (S3L) at (0,3) {S4};
\node[species-hidden] (S4L) at (0,2) {};
\node[species-inactive] (S5L) at (0,1) {S6};
\node[species-hidden] (S6L) at (0,0) {};

% --------------------
% REACCIONES
% --------------------
\node[reaction-inactive] (R0) at (2,5.5) {R1};
\node[reaction-hidden] (R1) at (2,4.5) {};
\node[reaction-active] (R2) at (2,3.5) {R3};
\node[reaction-inactive] (R3) at (2,2.5) {R4};
\node[reaction-inactive] (R4) at (2,1.5) {R5};
\node[reaction-hidden] (R5) at (2,0.5) {};

% --------------------
% NODOS DERECHA
% --------------------
\node[species-out-oldactive] (S0R) at (4,6) {S1};
\node[species-out-newactive] (S1R) at (4,5) {S2};
\node[species-hidden] (S2R) at (4,4) {};
\node[species-out-oldactive] (S3R) at (4,3) {S4};
\node[species-hidden] (S4R) at (4,2) {};
\node[species-hidden] (S5R) at (4,1) {};
\node[species-inactive] (S6R) at (4,0) {S7};

% --------------------
% ARISTAS ANTIGUAS (primero)
% --------------------
\draw[edge-in-inactive] (S3L) -- (R0);
\draw[edge-out-inactive]  (R0) -- (S1R);
\draw[edge-out-inactive]  (R0) -- (S0R);

\draw[edge-in-inactive] (S2L) -- (R3);
\draw[edge-in-inactive] (S0L) -- (R3);
\draw[edge-out-inactive]  (R3) -- (S0R);

\draw[edge-in-inactive] (S3L) -- (R4);
\draw[edge-in-inactive] (S0L) -- (R4);
\draw[edge-out-inactive]  (R4) -- (S6R);
\draw[edge-out-inactive]  (R4) -- (S3R);

% --------------------
% ARISTAS NUEVAS (después, por encima)
% --------------------
\draw[edge-in-active] (S1L) -- (R2);
\draw[edge-in-active] (S5L) -- (R2);
\draw[edge-out-active]  (R2) -- (S0R);
\draw[edge-out-active]  (R2) -- (S6R);

\end{tikzpicture}}
\caption{$t=3$}
\end{subfigure}
}
\caption{Solution of \eqref{gemd} for $3$ periods.}
\label{f:mak}
\end{figure}

The comparison between Figures~\ref{f:nomak} and \ref{f:mak} highlights the structural and temporal impact of \ref{gemd} versus \ref{geme}  on the example presented in the previous Figure~\ref{f:toy_example_t3} in Section \ref{sec:1}. 

Under \ref{geme} (Figure~\ref{f:nomak}), the expansion follows a more fragmented and flexible pattern: hyperarcs are activated based on immediate feasibility, leading to dispersed flows and limited coordination across periods. Different hyperarcs are sequentially selected without necessarily reinforcing a common pathway, resulting in weaker propagation effects but greater adaptability to local conditions.

In contrast, \ref{gemd} (Figure~\ref{f:mak}) induces a more coordinated behavior. The presence of synergistic constraints promotes the activation of hyperarcs whose input nodes can be jointly sustained, leading to coherent pathways (e.g., the persistence of hyperarcs $R1$, and $R5$ across periods), stronger flow concentration, and amplification effects across downstream nodes. This generates higher consistency over time and a more effective exploitation of complementarities.

Importantly, both approaches provide valuable modeling perspectives. GEM-E offers a more flexible and computationally tractable framework, as it avoids part of the nonlinearities induced by the synergistic interactions. Moreover, in many practical settings the parameters governing interaction intensities (e.g., the rates $\kappa$) may be unavailable, difficult to estimate, or lack a clear interpretation. In such cases, GEM-E constitutes a robust alternative that captures structural feasibility without requiring detailed kinetic information.

From a managerial perspective, the choice between both models reflects a trade-off between flexibility and coordination. While \ref{geme} supports adaptive and data-light decision-making, \ref{gemd} enables the identification of structured growth trajectories driven by complementarities. Together, they provide complementary tools to analyze and design expansion processes in complex systems.

The previous development has been presented under the simplifying assumption that no reversible interactions are considered in the hypergraph $\mathcal{H}=(\mathcal{N}, \mathcal{A})$. Informally, a hyperarc is said to be reversible if there exists another hyperarc that represents the same transformation but in the opposite direction, exchanging inputs and outputs, i.e., there are no hyperarcs $a, a' \in \A$ such that $S_a=T_{a'}$ and $S_{a'}=T_a$. While neglecting such interactions allows for a cleaner exposition, reversibility naturally arises in many applications. In particular, in chemical reaction networks, many reactions can proceed in both directions depending on the system conditions, leading to competing forward and backward processes. 

The framework can be extended to incorporate these reversible interactions without altering its fundamental structure. In this setting, reverse hyperarcs are treated as separate entities that contribute with opposite flow, capturing the net effect of bidirectional transformations. We assume that the underlying synergistic laws are not affected by the presence of these reverse interactions, and that their effect can be modeled through an appropriate decomposition of flows into forward and reverse components, together with suitable adjustments in the dynamic balance constraints. For completeness, the full mathematical treatment of this extension is provided in  \ref{app:reverse}.

% \red{A hyperarc $a' \in \A$ is said to be the reverse hyperarc of $a$ if $S_{a'} = T_a$ and $T_{a'} = S_a$. Let $r(\A) \subseteq \A$ denote the set of reversible hyperarcs. For simplicity, we assume $r(\A)=\emptyset$, although the general case is discussed in \ref{app:reverse}.}

\section{Computational Experiments}\label{sec:4}

This section presents a series of computational experiments designed to evaluate the performance, interpretability, and applicability of the proposed GEM formulations. The experiments aim to illustrate how the models capture the temporal evolution of cooperative amplification processes across multiple domains. We first analyze small illustrative instances to visualize the incremental activation of nodes and hyperarcs and to validate the logical consistency of the model constraints. We then examine larger synthetic and application-inspired test cases to assess computational scalability, the impact of nonlinear flow interactions, and the influence of equity-oriented design criteria on the resulting hypergraph structures.

To ensure reproducibility and structural diversity, the experimental data is generated using the \texttt{SMGen} instance generator~\citep{Riva2022SMGen}, available at \url{https://gitlab.com/simoriva/smgen}. This tool produces realistic synthetic chemical reaction networks (and, equivalently, directed hypergraphs) under controlled topological and kinetic (synergistic rate) parameters, offering a flexible and systematic means to benchmark model behavior across varying network sizes and reaction densities. In particular, we consider stoichiometric matrices of size $n \times n$, with $n \in \{25, 50, 75, 100\}$, and control the density of the instances through an upper bound $d \in \{2,3,5\}$ on the number of species involved in each reaction, both on the input and output sides of the corresponding hyperarcs. \texttt{SMGen} also provides randomly generated synergistic coefficients $\kappa$ (kinetic constants) under different probabilistic distributions. We choose $\kappa$ to be independently sampled from a uniform distribution on $[0,1]$. For each combination of parameters, five random instances are generated.

We evaluate both frameworks, \eqref{geme} and \eqref{gemd}, considering different levels of expansion and temporal depth. Specifically, we set the number of activated components to $q \in \{1,2\}$ and the number of periods to $T \in \{1,2,3,4,5\}$. 
All experiments were run on the Linux-based Supercomputer \texttt{albaicin}, composed of 9520 Intel\textsuperscript{\textregistered} Xeon\textsuperscript{\textregistered} Cores at 2.7 GHz and 35 TB RAM total, as 170 nodes interconnected with Infiniband HDR200 network. We use 16 threads with 8 gigabytes of RAM for our experiments. The models were coded in Python 3.8.2 with Gurobi 12.1 as optimization solver, and a time limit of 2 hours is fixed for all the instances.

Table~\ref{tab:runtime_geme} reports the median  computational times required to solve instances of the \ref{geme} formulation across different problem sizes and density levels for a single period and with $q = 1$. Overall, the results indicate that the proposed mixed-integer formulation can be solved very efficiently for the considered instance sizes, with runtime consistently below a fraction of a second in most cases. This suggests that the structural constraints imposed by the model, including nested activation and flow balance, do not introduce significant computational demand at moderate scales.

Regarding scalability, no clear monotonic increase in runtime is observed as the number of species grows, which highlights the robustness of the formulation with respect to problem size. The effect of density is also limited, although a noticeable increase appears in the instance with 75 species and density 5, suggesting that higher interaction complexity may lead to occasional computational spikes. Nevertheless, even in this case, solution times remain within a very small range from a practical perspective. 

These results provide evidence that \ref{geme} can be effectively applied to medium-scale instances, supporting its potential use in real-world applications where timely decision-making is required.
\begin{table}[h]
\centering
\begin{tabular}{cccc}
\toprule
\multirow{2}{*}{Size} & \multicolumn{3}{c}{Density} \\
\cmidrule(lr){2-4}
 & 2 & 3 & 5 \\
\midrule
$25 \times 25$  & 0.0211 & 0.0184 & 0.0181 \\
$50 \times 50$ & 0.0193 & 0.0190 & 0.0211 \\
$75 \times 75$  & 0.0187 & 0.0108 & 0.1167 \\
$100 \times 100$ & 0.0212 & 0.0252 & 0.0256 \\
\bottomrule
\end{tabular}
\caption{Median runtime (in seconds) of the \ref{geme} formulation for different problem sizes and densities for a single period and with $q = 1$.}
\label{tab:runtime_geme}
\end{table}

Table~\ref{tab:geme_q1} reports the median runtime of the \ref{geme} formulation for $q=1$ across different matrix sizes, densities, and time horizons, excluding the presented single period. The results show that all instances are solved extremely efficiently, with runtime remaining below a few seconds even in the largest and most demanding configurations.

From a scalability perspective, the number of periods has the most noticeable impact on computational time. While instances with $t \leq 3$ are solved almost instantaneously across all sizes and densities, larger time horizons ($t=4,5$) lead to a moderate increase in runtime, particularly for higher densities. This reflects the growth in the number of activation decisions and intertemporal constraints as the planning horizon expands.

Density also plays a significant role, with higher values ($d=5$) consistently leading to larger runtime. This is expected, as denser hypergraphs induce more complex interaction structures and increase the combinatorial difficulty of identifying self-amplifying substructures. In contrast, the effect of the number of species is more gradual, suggesting that the formulation scales well with problem size when interaction density is moderate.

Overall, the results indicate that \ref{geme} remains computationally tractable across all tested configurations, supporting its applicability to medium-scale problems involving dynamic hypergraph structures.
\begin{table}[htbp]
\centering
\begin{tabular}{cccccc}
\toprule
Periods & Size & $d=2$ & $d=3$ & $d=5$ \\
\midrule
\multirow{4}{*}{2}
& $25\times25$  & 0.0214 & 0.0246 & 0.0227 \\
& $50\times50$  & 0.0226 & 0.0315 & 0.0306 \\
& $75\times75$  & 0.0230 & 0.1534 & 0.2305 \\
& $100\times100$& 0.0787 & 0.0342 & 0.0433 \\
\midrule
\multirow{4}{*}{3}
& $25\times25$  & 0.0266 & 0.0313 & 0.0295 \\
& $50\times50$  & 0.0257 & 0.0382 & 0.0608 \\
& $75\times75$  & 0.0283 & 0.0509 & 0.7102 \\
& $100\times100$& 0.1579 & 0.0747 & 0.0847 \\
\midrule
\multirow{4}{*}{4}
& $25\times25$  & 0.0741 & 0.1524 & 0.2261 \\
& $50\times50$  & 0.0294 & 0.0925 & 0.3600 \\
& $75\times75$  & 0.0396 & 0.3831 & 0.7078 \\
& $100\times100$& 0.0716 & 0.4858 & 0.9735 \\
\midrule

\multirow{4}{*}{5}
& $25\times25$  & 0.0420 & 0.2476 & 0.0968 \\
& $50\times50$  & 0.0337 & 0.0777 & 0.6433 \\
& $75\times75$  & 0.0489 & 0.1299 & 1.1520 \\
& $100\times100$& 0.0992 & 1.3474 & 1.1035 \\
\bottomrule
\end{tabular}
\caption{Median runtime (in seconds) of \ref{geme} for $q=1$.}
\label{tab:geme_q1}
\end{table}

The results for $q=2$, in Table \ref{tab:geme_q2}, exhibit a clear increase in computational effort compared to the case $q=1$, reflecting the additional combinatorial complexity induced by enforcing stronger balance requirements across periods. While small instances remain solvable within negligible time, the runtime grows more noticeably with both the number of periods and the density parameter. In particular, high-density instances ($d=5$) and larger sizes ($n \geq 75$) lead to a substantial increase in solution times, reaching several seconds for longer horizons. This behavior is consistent with the tighter coupling between decisions across periods introduced by the $q=2$ objective, which reduces flexibility and increases the difficulty of identifying feasible growth trajectories. Nevertheless, all instances remain tractable within very short computational times, confirming that the proposed formulation can efficiently handle moderate-size hypergraphs even under stronger equity or balance requirements. These results highlight a trade-off between modeling richness and computational effort, which is relevant for practitioners when selecting the appropriate level of temporal coordination in applications.
\begin{table}[htbp]
\centering
\begin{tabular}{cccccc}
\toprule
Periods & Size & $d=2$ & $d=3$ & $d=5$ \\
\midrule

\multirow{4}{*}{2}
& $25\times25$  & 0.0355 & 0.0186 & 0.0118 \\
& $50\times50$  & 0.0217 & 0.0143 & 0.0304 \\
& $75\times75$  & 0.0201 & 0.0182 & 0.1659 \\
& $100\times100$& 0.0220 & 0.0238 & 0.0308 \\
\midrule

\multirow{4}{*}{3}
& $25\times25$  & 0.0201 & 0.1142 & 0.0190 \\
& $50\times50$  & 0.0211 & 0.1208 & 0.0773 \\
& $75\times75$  & 0.0949 & 0.1724 & 0.6046 \\
& $100\times100$& 0.1243 & 0.5068 & 0.2203 \\
\midrule

\multirow{4}{*}{4}
& $25\times25$  & 0.0575 & 0.1612 & 0.3358 \\
& $50\times50$  & 0.0477 & 0.2700 & 0.6366 \\
& $75\times75$  & 0.0422 & 0.1834 & 0.9166 \\
& $100\times100$& 0.0609 & 0.5926 & 0.1394 \\
\midrule

\multirow{4}{*}{5}
& $25\times25$  & 0.0917 & 0.2583 & 0.2879 \\
& $50\times50$  & 0.0595 & 0.5977 & 2.9870 \\
& $75\times75$  & 0.1352 & 3.1963 & 5.6468 \\
& $100\times100$& 0.1862 & 1.0174 & 4.6984 \\
\bottomrule
\end{tabular}
\caption{Median runtime (in seconds) of \ref{geme} for $q=2$.}
\label{tab:geme_q2}
\end{table}
Compared to the $q=1$ setting, the increase in runtime is more pronounced for dense and long-horizon instances, indicating that the choice of $q$ plays a critical role in the scalability of the model.

The results for $n=100$ (Table \ref{tab:geme_n100_q}) provide additional insight into the impact of the parameter $q$ and the planning horizon on computational performance. As expected, the runtime increases with the number of periods, reflecting the growing temporal coupling and the expansion of the feasible space. For $q=1$, the growth is smooth and remains below one second even for five periods. In contrast, the case $q=2$ exhibits a more irregular but generally higher computational demand, particularly for longer horizons, where the runtime slightly exceeds one second. This behavior is consistent with the more restrictive nature of the $q=2$ objective, which enforces stronger balance conditions across periods and reduces flexibility in the selection of activation patterns. Interestingly, for intermediate horizons, the performance does not deteriorate monotonically, suggesting that the interaction between temporal depth and combinatorial structure may create instances with varying difficulty. Overall, the results confirm that the model remains computationally tractable for realistic sizes while highlighting the influence of $q$ as a key parameter controlling both modeling richness and solution complexity.
\begin{table}[htbp]
\centering
\begin{tabular}{c|ccccc}
\toprule
$q$ $\backslash$ Periods & 1 & 2 & 3 & 4 & 5 \\
\midrule
1 & 0.0252 & 0.0395 & 0.0847 & 0.3824 & 0.8701 \\
2 & N/A     & 0.0276 & 0.2031 & 0.1394 & 1.0174 \\
\bottomrule
\end{tabular}
\caption{Median runtime (in seconds) of \ref{geme} for $n=100$ under different values of $q$ and time horizons.}
\label{tab:geme_n100_q}
\end{table}

Table~\ref{tab:gemd_avg} reports the median runtime and percentage of unsolved instances for \ref{gemd} across different problem sizes, densities, and time horizons. The results clearly illustrate the significant impact of both density and temporal depth on computational performance. For single-period instances, the model remains tractable across all sizes, although runtimes increase notably with density, especially for $d=5$. However, when extending the horizon to two periods, a sharp escalation in computational difficulty is observed. In particular, for medium and high densities ($d=3$ and $d=5$), median runtime frequently approach the imposed time limit, and the proportion of unsolved instances rises dramatically, reaching up to $100\%$ in several configurations.

This behavior reveals a strong interaction effect between density and temporal coupling: while each dimension alone remains manageable, their combination leads to a substantial increase in complexity. Larger instances exacerbate this effect, particularly under high-density settings, where even moderate increases in size or horizon result in a loss of solvability. These findings highlight a clear trade-off between modeling richness and computational tractability. From a practical standpoint, they suggest that \ref{gemd} is best suited for moderate-size or carefully selected instances, whereas high-density and multi-period configurations may require tailored solution strategies or approximations to remain computationally viable.
\begin{table}[h]
\centering
\begin{tabular}{cc|cc|cc|cc}
\toprule
& & \multicolumn{2}{c|}{$d=2$} & \multicolumn{2}{c|}{$d=3$} & \multicolumn{2}{c}{$d=5$} \\
\cline{3-8}
Size & $T$ 
& Time & UnS(\%) 
& Time & UnS(\%) 
& Time & UnS(\%) \\
\midrule

\multirow{2}{*}{$25\times25$}
& 1 & 0.0726 & 0   & 0.2775 & 0   & 3.0862 & 0 \\
& 2 & 0.2931 & 0   & 2880.4764 & 40 & 6084.0247 & 80 \\
\midrule

\multirow{2}{*}{$50\times50$}
& 1 & 0.0601 & 0   & 1.5812 & 0   & 239.1730 & 0 \\
& 2 & 0.2481 & 0   & 6107.85 & 80 & \texttt{TL}  & 100 \\
\midrule

\multirow{2}{*}{$75\times75$}
& 1 & 0.1119 & 0   & 1.8240 & 0   & 1564.1453 & 20 \\
& 2 & 1800.15 & 25 & 5763.72 & 80 & \texttt{TL}  & 100 \\
\midrule

\multirow{2}{*}{$100\times100$}
& 1 & 0.1087 & 0   & 1456.1442 & 20 & 3603.76 & 50 \\
& 2 & 30.7507 & 0   & \texttt{TL}  & 100 & \texttt{TL} & 100 \\
\bottomrule
\end{tabular}
\caption{Median runtime (seconds) and percentage of unsolved instances for \ref{gemd} across different sizes, densities ($d$), and time horizons.}
\label{tab:gemd_avg}
\end{table}
Overall, Table \ref{tab:gemd_avg} provides strong evidence of a phase-transition behavior in problem difficulty, driven by the joint effect of density and temporal expansion.

\section{Case Study: Economic Production Networks}\label{sec:5}

To illustrate the applicability of the proposed framework, we consider an economic production network derived from the Input--Output Industry Economic Accounts published by the \textit{U.S.\ Bureau of Economic Analysis (BEA)}\footnote{Datasets are publicly available at \url{https://apps.bea.gov}.}. We focus on the benchmark year 2017, which provides the most detailed publicly available estimates.

The BEA accounts are released at four levels of aggregation: \emph{sector} (21 industries), \emph{summary} (71 industries), \emph{underlying summary} (138 industries), and \emph{detail} (402 industries). These levels follow a hierarchical structure consistent with the 2017 North American Industry Classification System (NAICS), where increasingly disaggregated levels capture finer technological and input--output relationships. This multi-resolution structure enables the construction of alternative network representations, allowing us to assess the robustness of our methodology across different levels of economic granularity.

The dataset is based on the \emph{Use} and \emph{Make} tables at producer prices. The Use table specifies how commodities are consumed as intermediate inputs by industries, while the Make table records the commodities produced by each industry. These tables naturally define the matrices $\mathbb{S}$ and $\mathbb{T}$ (and hence $\mathbb{Q}=\mathbb{T}-\mathbb{S}$), which constitute the fundamental inputs for constructing the corresponding production hypergraph.

Within this representation, \emph{nodes} correspond to commodities (goods and services), while \emph{hyperarcs} encode production processes. Each hyperarc links a set of input commodities to a set of output commodities, capturing the transformation mechanisms carried out by industries. This modeling paradigm allows us to naturally represent higher-order interactions, where multiple inputs are jointly transformed into multiple outputs.

This structure is particularly suitable for capturing complex economic interdependencies. In economic terms, self-amplifying subnetworks correspond to mutually reinforcing sectors in which the production of certain commodities sustains and expands the activities that generate them. For instance, manufacturing outputs are used by trade and transportation sectors, which in turn provide essential services that support further production. These feedback loops form reinforcing cycles that sustain and amplify economic activity. Within our framework, such phenomena are precisely captured by expanding self-amplifying hypergraphs.

We could construct four alternative datasets corresponding to the different four different levels of aggregation. For visualization purposes, we only use here the most aggregated setting, 21 sectors, which provides a coarse representation of broad economic activities (e.g., agriculture, manufacturing, services). In Table \ref{t:icons}, each sector is represented using an intuitive icon to facilitate interpretation.

\begin{table}[ht]
\small
\centering
\renewcommand{\arraystretch}{1.2}
\begin{tabular}{p{0.05\textwidth} p{0.38\textwidth} p{0.05\textwidth} p{0.38\textwidth}}

\multicolumn{2}{c}{\textit{Primary activities}} & \multicolumn{2}{c}{\textit{Industrial production}} \\

\faIcon{seedling} & Agriculture, Forestry, Fishing, and Hunting &
\faIcon{bolt} & Utilities \\

\faIcon{hard-hat} & Mining &
\faIcon{tools} & Construction \\

& &
\faIcon{industry} & Durable Goods Manufacturing \\

& &
\faIcon{box} & Nondurable Goods Manufacturing \\[0.3em]

\multicolumn{2}{c}{\textit{Trade and transportation}} & \multicolumn{2}{c}{\textit{Business and financial services}} \\

\faIcon{truck-loading} & Wholesale Trade &
\faIcon{satellite} & Information \\

\faIcon{shopping-cart} & Retail Trade &
\faIcon{coins} & Finance and Insurance \\

\faIcon{truck} & Transportation and Warehousing &
\faIcon{home} & Real Estate and Rental and Leasing \\

& &
\faIcon{brain} & Professional and Technical Services \\

& &
\faIcon{building} & Management of Companies and Enterprises \\

& &
\faIcon{trash} & Administrative and Waste Services \\[0.3em]

\multicolumn{2}{c}{\textit{Social and personal services}} & \multicolumn{2}{c}{\textit{Public sector}} \\

\faIcon{graduation-cap} & Educational Services &
\faIcon{landmark} & Government \\

\faIcon{heartbeat} & Health Care and Social Assistance & & \\

\faIcon{theater-masks} & Arts, Entertainment, and Recreation &
\multicolumn{2}{c}{\textit{Others}} \\

\faIcon{utensils} & Accommodation and Food Services &
\faIcon{recycle} & \textit{Scrap, Used, and Secondhand Goods} \\

\faIcon{cogs} & Other Services (except Government) &
\faIcon{globe} & \textit{Noncomparable Imports and Rest-of-the-World Adjustment} \\

\end{tabular}
\caption{Icon labeling for the 21 sectors in the BEA dataset.}\label{t:icons}
\end{table}

Although the BEA sector classification includes 21 industries, these additional categories are necessary for a consistent representation of economic flows. Specifically, \faIcon{recycle} represents recycled and secondhand goods re-entering production processes, while \faIcon{globe} captures balancing adjustments associated with international trade.

From the perspective of our hypergraph formulation, these categories do not behave as standard commodity nodes. Instead, they act as transformation mechanisms redistributing flows across the network. Consequently, we model them as \emph{hyperarcs} rather than nodes, ensuring a faithful representation of the underlying accounting structure.

\begin{figure}[htbp]
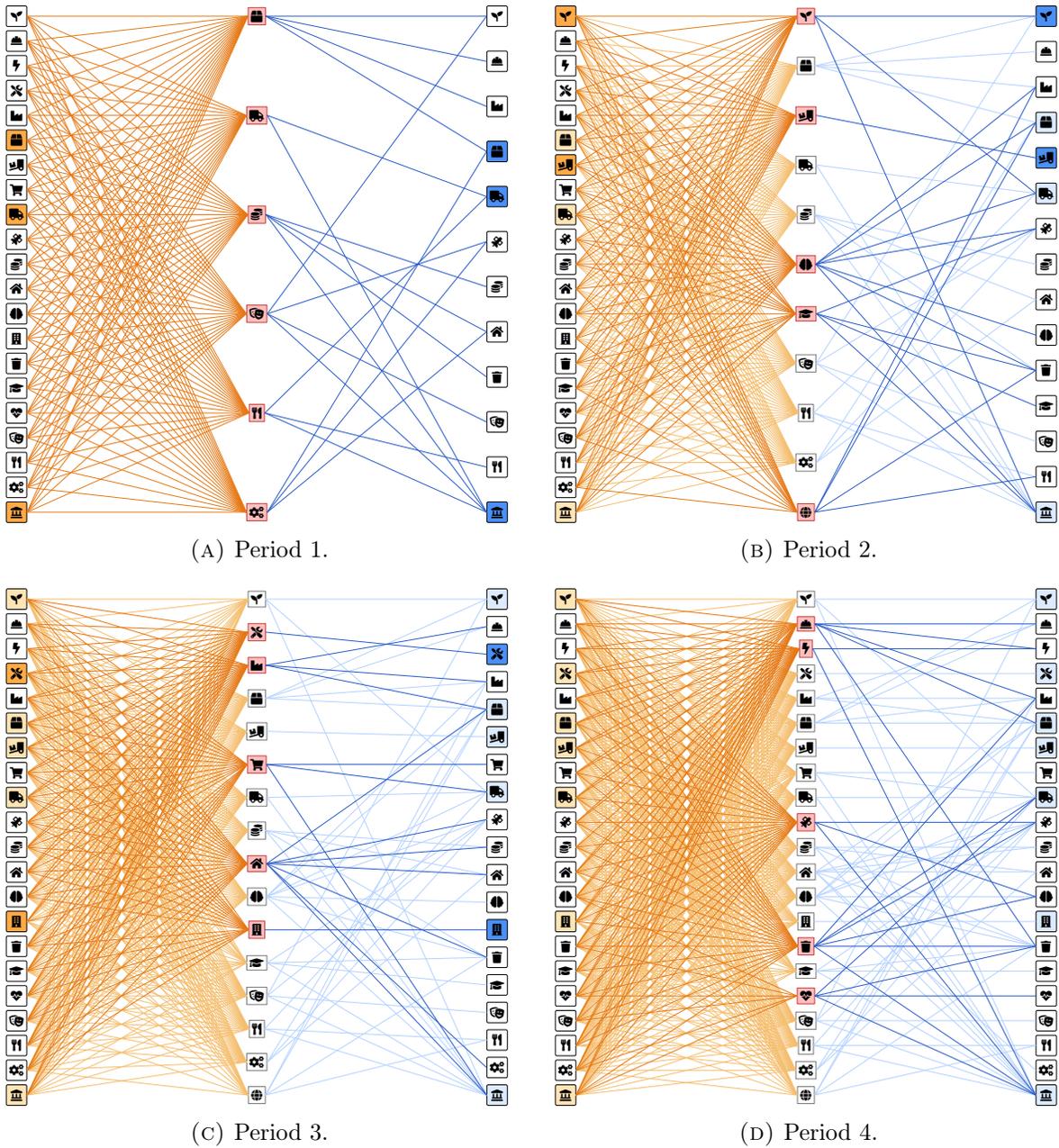

    \centering

    \begin{subfigure}[t]{0.48\textwidth}
        \centering
        \resizebox{\linewidth}{!}{\input{cs_BEA_0_4_1_1_001_t1}}
        \caption{Period 1.}
        \label{f:bea_p1}
    \end{subfigure}
    \hfill
    \begin{subfigure}[t]{0.48\textwidth}
        \centering
        \resizebox{\linewidth}{!}{\input{cs_BEA_0_4_1_1_001_t2}}
        \caption{Period 2.}
        \label{f:bea_p2}
    \end{subfigure}

    \vspace{0.8em}

    \begin{subfigure}[t]{0.48\textwidth}
        \centering
        \resizebox{\linewidth}{!}{\input{cs_BEA_0_4_1_1_001_t3}}
        \caption{Period 3.}
        \label{f:bea_p3}
    \end{subfigure}
    \hfill
    \begin{subfigure}[t]{0.48\textwidth}
        \centering
        \resizebox{\linewidth}{!}{\input{cs_BEA_0_4_1_1_001_t4}}
        \caption{Period 4.}
        \label{f:bea_p4}
    \end{subfigure}

    \caption{Tripartite graph representation of the solutions obtained for the four periods in the BEA case study. Each subfigure shows the activated sectors and reactions at the corresponding period, allowing the progressive construction of the autocatalytic production network to be visualized.}
    \label{f:cs_BEA}
\end{figure}

The corresponding solutions are depicted in Figure~\ref{f:cs_BEA} using a tripartite graph, where the evolution of activated production processes and sectoral interdependencies can be tracked across periods. As in the illustrative example presented earlier, input commodities are placed on the left, production processes (hyperarcs) are represented in the center, and output commodities are shown on the right.

Colored nodes denote commodities that belong to the \emph{self-amplifying production structure}, whereas lightly colored nodes indicate commodities that were already active in previous periods. Similarly, production processes highlighted in red correspond to newly activated processes, while those shown in white represent processes that had already been activated in earlier periods.

We then run our model over four consecutive periods to analyze the progressive construction of a \emph{self-amplifying production network} under the non-synergetic formulation. Although intensity parameters could in principle be derived from the weights of the input and output matrices, we deliberately avoid this specification in the present case study. The economic production network considered here already constitutes a challenging instance due to its high density and the complexity of its intersectoral dependencies, as discussed below.

In the context of this dataset, a self-amplifying production subnetwork refers to a subset of commodities and production processes such that the endogenous interactions among them generate a net expansion of the commodities within the subset. In other words, the selected production processes transform available inputs into outputs in a way that sustains and reinforces the same set of commodities over time without requiring additional external support beyond the initial activation. From an economic perspective, these subnetworks capture groups of sectors whose mutual interdependencies enable persistent growth through internal feedback mechanisms.

The first period identifies the set of commodities and production processes that initiate the expansion of the system. Since the data correspond to the year 2017, this initial configuration should be interpreted as the state inherited from the end of 2016. Therefore, the activated commodities in this period do not arise from an empty system but rather represent those sectors that were already sufficiently established to act as drivers of further production activity. In the obtained results, the initial self-amplifying commodities are Nondurable Goods Manufacturing (\faIcon{box}), Transportation and Warehousing (\faIcon{truck}), and Government (\faIcon{landmark}). These sectors constitute the first productive core from which the structure begins to expand. At the same time, the activated production processes already involve a broader set of activities, including Nondurable Goods Manufacturing (\faIcon{box}), Transportation and Warehousing (\faIcon{truck}), Finance and Insurance (\faIcon{coins}), Arts, Entertainment, and Recreation (\faIcon{theater-masks}), Accommodation and Food Services (\faIcon{utensils}), and Other Services (\faIcon{cogs}). This shows that, even at the initial stage, the economy contains a richer set of operative transformations than the limited set of commodities that can already be sustained in a self-amplifying way. Economically, this means that the inherited configuration at the end of 2016 already supports a nontrivial productive base, but only a few sectors display positive internal reinforcement.

As the system evolves over the subsequent periods, the self-amplifying production structure expands through the activation of additional production processes and the incorporation of new commodities. In the second period, the set of self-amplifying commodities expands to include Agriculture, Forestry, Fishing, and Hunting (\faIcon{seedling}) and Wholesale Trade (\faIcon{truck-loading}), together with the previously active commodities. This broadening of the active commodity base is accompanied by the activation of new production processes in Agriculture (\faIcon{seedling}), Wholesale Trade (\faIcon{truck-loading}), Professional and Technical Services (\faIcon{brain}), Educational Services (\faIcon{graduation-cap}), and the external adjustment category represented by Noncomparable Imports and Rest of the World Adjustment (\faIcon{globe}). The results suggest that once the initial core is in place, the productive structure begins to absorb upstream supply activities and business support functions that reinforce the circulation of goods and services across the network. Hence, the second period captures the transition from a narrow productive core to a more articulated structure in which primary activities, intermediate trade, and specialized services start to play a reinforcing role.

This expansion continues in the third period, where two additional commodities become self-amplifying, namely Construction (\faIcon{tools}) and Management of Companies and Enterprises (\faIcon{building}). The newly activated production processes in this period are Construction (\faIcon{tools}), Durable Goods Manufacturing (\faIcon{industry}), Retail Trade (\faIcon{shopping-cart}), Real Estate and Rental and Leasing (\faIcon{home}), and Management of Companies and Enterprises (\faIcon{building}). This pattern is economically meaningful, since it indicates a progressive densification of the system from core logistics and public demand towards sectors associated with infrastructure, coordination, commercialization, and asset management. In other words, the self-amplifying structure evolves from a relatively narrow base into a more mature production network in which goods movement, intermediation, and organizational capabilities are increasingly integrated.

In the fourth period, no new commodities are incorporated into the self-amplifying structure. However, additional production processes are activated, namely Mining (\faIcon{hard-hat}), Utilities (\faIcon{bolt}), Information (\faIcon{satellite}), Administrative and Waste Services (\faIcon{trash}), and Health Care and Social Assistance (\faIcon{heartbeat}). This confirms that the system has reached a saturation point in terms of the set of commodities that can sustain positive internal reinforcement. Further expansion is therefore achieved through the activation of new production processes that enrich the transformation possibilities among the already active commodities rather than through the inclusion of additional commodities. From an economic viewpoint, this indicates that the productive base has already been identified by the third period, and that subsequent growth depends on deepening the connectivity, service intensity, and functional sophistication of the existing structure. The late activation of processes such as Utilities, Information, and Health Care is also revealing, since these sectors appear as important supporting mechanisms for the productive core, even though they do not themselves become self-amplifying commodities within the time horizon considered.

Finally, the commodities that appear in the graph as inputs or outputs but are not activated as self-amplifying can, in general, be interpreted in three possible ways. Some of them may behave as \textit{food} commodities, that is, sectors that are required as inputs but are not reproduced within the selected structure. Others may behave as \textit{waste} commodities, namely sectors that are generated as outputs but do not feed back into the productive core. A third possibility corresponds to commodities that appear both as inputs and outputs but still exhibit negative net production, so that they cannot be sustained internally despite participating in several processes. 

Depending on the instance under analysis, one or more of these categories may arise, or some of them may be absent altogether. This distinction is useful because it clarifies which sectors constitute the endogenous engine of growth, which ones may act as external resources, and which ones remain structurally dependent on the rest of the economy.

To complement the graphical analysis, Table~\ref{tab:classification} summarizes the classification of the commodities according to their role within the identified structure. In particular, we distinguish between self-amplifying commodities, which constitute the endogenous core of the system, and those commodities that, although present in the network, exhibit negative net production and therefore cannot be sustained internally.

The set of self-amplifying commodities captures the sectors that are able to reproduce themselves through the activated production processes, thus forming the backbone of the productive system. In contrast, the commodities with negative net production correspond to sectors that participate in multiple transformations but whose total output remains insufficient to compensate for their consumption within the selected subnetwork.

From an economic perspective, this behavior can be explained by the functional role of these sectors. Many of the negative net production commodities are associated with supporting, intermediate, or service-oriented activities, such as financial services, information, retail, or health care. These sectors are essential for enabling and coordinating production, distribution, and consumption, but they are not primary drivers of material or value reproduction within the restricted structure identified by the model. As a consequence, they depend on external inputs or on interactions outside the selected subnetwork to remain viable.

This observation highlights an important structural feature of the system: the endogenous growth is sustained by a relatively small subset of sectors, while a broader set of complementary activities plays a facilitating role. The distinction between these two groups provides valuable managerial insight, as it identifies the industries that act as engines of expansion and those that, although indispensable, require external support or broader system integration to operate sustainably.

\begin{table}[h]
\begin{center}
\renewcommand{\arraystretch}{1.15}
\begin{tabular}{ll}
\textbf{Type} & \textbf{Detected commodities} \\ \hline
{ Self-amplifying}
& \faIcon{seedling} \;
  \faIcon{tools} \;
  \faIcon{box} \;
  \faIcon{truck-loading}\;
  \faIcon{truck} \;
  \faIcon{building} \;
  \faIcon{landmark} \\\hline
{ Negative net production}
& \faIcon{hard-hat} \;
  \faIcon{bolt} \;\;\,
  \faIcon{industry} \;
  \faIcon{shopping-cart} \;
  \faIcon{satellite} \;
  \faIcon{coins} \,
  \faIcon{home} \\
 & \faIcon{brain} \;
  \faIcon{trash} \;
  \faIcon{graduation-cap} \;
  \faIcon{heartbeat} \,
  \faIcon{theater-masks} \,
  \faIcon{utensils} \;
  \faIcon{cogs} 
\end{tabular}
\end{center}
\caption{Classification of commodities after the last period.\label{tab:classification}}
\end{table}

Overall, in view of this case study, the proposed framework provides a practical and interpretable tool for practitioners and decision-makers concerned with strategic investment and economic planning. By explicitly identifying the set of self-amplifying commodities, the model reveals the fundamental sectors that act as endogenous engines of growth, thereby offering a principled basis for prioritizing investments and resource allocation. At the same time, the multiperiod representation delivers a dynamic view of how these cores emerge and expand, allowing policymakers to anticipate which sectors are likely to become critical in subsequent stages of development. The graphical visualization of the structure across periods further enhances interpretability, as it provides an intuitive and transparent depiction of the evolving production network and the role played by each commodity. This combination of mathematical rigor and visual insight facilitates the translation of complex interdependencies into actionable knowledge, making the approach particularly valuable for supporting evidence-based decisions in economic policy, industrial strategy, and infrastructure planning.

\section{Conclusions and Further Research}\label{sec:6}

This paper introduces a novel optimization-based framework for the analysis of directed multiset hypergraphs, with a particular emphasis on the identification of self-amplifying structures. By adopting a multiperiod perspective, the proposed approach captures the progressive activation of nodes and hyperarcs, providing a dynamic and interpretable representation of how endogenous production cores emerge and expand over time. This temporal and nested viewpoint moves beyond static network analysis and offers a richer description of growth mechanisms in complex systems.

From an Operations Research perspective, the main contribution lies in the integration of higher-order network modeling with advanced optimization techniques. We formalize the detection of amplification structures as optimization problems and show that they can be addressed through mixed-integer linear formulations in the non-synergetic setting, as well as mixed-integer nonlinear models when incorporating synergistic flow laws that mimic mass-action kinetics. This dual modeling capability enables practitioners to balance tractability and modeling fidelity, providing a unified framework that connects combinatorial structure selection with nonlinear flow dynamics. The resulting models allow for the systematic extraction of economically meaningful substructures, while at the same time establishing a direct link with self-amplifying networks studied in chemistry, particularly in the context of autocatalytic systems and their role in theories on the Origin of Life~\citep{hordijk2018autocatalytic}. In this sense, the proposed framework not only advances OR methodology but also provides a quantitative tool for analyzing the emergence and sustainability of autocatalytic structures in complex reaction networks.

From a managerial and economic standpoint, the framework provides actionable insights for decision-makers concerned with resource allocation, industrial policy, and strategic planning. In particular, the identification of self-amplifying commodities highlights the sectors that constitute the endogenous engine of growth, thereby offering a principled basis for prioritizing investments. At the same time, the multiperiod structure reveals how these cores develop over time, allowing policymakers to anticipate future critical sectors and to design phased intervention strategies. The associated visualizations further enhance interpretability by providing an intuitive representation of the evolving production network, facilitating the translation of complex interdependencies into transparent and actionable knowledge.

The results of this work open several promising directions for future research. First, incorporating uncertainty into the hypergraph framework is a natural extension. In many real-world applications, input--output coefficients, technological efficiencies, and resource availability are subject to variability. This motivates the development of stochastic and distributionally robust optimization models capable of characterizing amplification structures under uncertainty and supporting robust decision-making.

Second, the interpretation of large-scale multiperiod hypergraphs remains a challenge. The development of advanced visualization tools, supported by optimization-based layout design, could significantly enhance the understanding of nested structures, feedback mechanisms, and temporal dependencies, particularly in high-dimensional settings.

Third, from a computational viewpoint, the proposed models may become demanding as the size and temporal depth of the hypergraph increase. This calls for the design of tailored decomposition strategies, such as Benders decomposition across periods, column generation over hyperarcs, or Lagrangian relaxation of coupling constraints. The development of efficient heuristics and approximation schemes will be essential for scaling the approach to large industrial applications.

Finally, several broader research avenues emerge from this framework. These include the study of adaptive and evolving hypergraphs, the integration of data-driven methods to estimate interaction parameters and activation dynamics, and the exploration of connections with growth theory in economics, autocatalytic sets in chemical reaction networks, and network design problems in logistics and infrastructure systems. Overall, the proposed framework provides a unifying and optimization-driven perspective for analyzing higher-order interactions, offering both methodological advances in Operations Research and practical tools for decision-making in complex environments.

\section*{Acknowledgements}

Authors thank the ``Servicio de Supercomputacion'' from Universidad de Granada (\url{https://supercomputacion.ugr.es}) for providing computing time on the albaicin supercomputer. 
The authors acknowledge financial support by  grants PID2020-114594GB-C21 and  
PID2024-156594NB-C21 funded by MICIU/AEI/10.13039/501100011033; FEDER+Junta de Andalucía projects C‐EXP‐139‐UGR23, and  the IMAG-Maria de Maeztu grant CEX2020-001105-M / AEI / 10.13039 / 501100011033. The third author was also supported by the Ph.D. Program in Mathematics at the Universidad de Granada.

\bibliographystyle{plainnat}
\bibliography{multiperiod_hypergraph}

\newpage
\appendix

\section{Reformulation of logarithmic constraints \eqref{ctr:12}}\label{app:1}

 Note that constraints~\eqref{ctr:12} are highly nonlinear. In what follows, we derive a suitable reformulation for the synergistic flow law conditions based on a logarithmic representation.

Note that for $z_a^t=1$, constraint~\eqref{ctr:12} is equivalent to
\begin{equation}
\log(f_a^t) = \log(\kappa_a) + \sum_{v\in\mathcal{N}: x_v^t>0} \mathbb{S}_{va} \log(x_v^t).\label{logf}
\end{equation}

To handle zero values, we introduce auxiliary variables:
$$
h_v^t =
\begin{cases}
\log(x_v^t) & \text{if node $v$ is active,}\\
0 & \text{otherwise,}
\end{cases}
\quad
g_a^t =
\begin{cases}
\log(f_a^t) & \text{if hyperarc $a$ is active,}\\
0 & \text{otherwise.}
\end{cases}
$$

These are enforced by:
\begin{align}
    h_v^t &= \log(x_v^t - \rho_v^t + 1), \label{ctr:13a}\\
    g_a^t &= \log(f_a^t - z_a^t + 1). \label{ctr:13b}
\end{align}

And then, \eqref{logf} is enforced by:
\begin{align}
g_a^t \geq \log(\kappa_a) + \sum_{v} \mathbb{S}_{va} h_v^t - \Delta_a (1-z_a^t), \label{ctr:14a}\\
g_a^t \leq \log(\kappa_a) + \sum_{v} \mathbb{S}_{va} h_v^t + \Delta_a (1-z_a^t). \label{ctr:14b}
\end{align}

Note also that constraints \eqref{ctr:13a}--\eqref{ctr:13b} can be approximated via SOS2 variables using breakpoints:
$$
L_v = \xi_{v,1} < \cdots < \xi_{v,K} = U_v,
$$
and convex combinations:
\begin{align*}
x_v^t - \rho_v^t + 1 = \sum_k \xi_{v,k} \lambda_{v,k}^t,\\
h_v^t = \sum_k \log(\xi_{v,k}) \lambda_{v,k}^t.
\end{align*}

The variables $\lambda_{v,k}^t$ satisfy $\sum_k \lambda_{v,k}^t = 1$, $\lambda_{v,k}^t \geq 0$, and are restricted to form a special ordered set of type 2 (SOS2), ensuring that at most two consecutive breakpoints are active. This structure guarantees that $x_v^t - \rho_v^t + 1$ is represented as a convex combination of adjacent breakpoints, yielding a piecewise linear interpolation of the logarithmic function.

Analogous constructions can be introduced for constraints \eqref{ctr:13b}, defining appropriate breakpoint sets for the flow variables $f_a^t$ and associated SOS2 variables. In this way, all logarithmic terms in the model can be consistently approximated within a unified linear framework.

This piecewise linear reformulation preserves the convexity of the logarithmic functions and provides a tight outer approximation of the original nonlinear constraints. The approximation accuracy can be controlled by the number and placement of breakpoints: a finer discretization improves fidelity at the cost of increased model size.

\section{Extension to Reversible Hyperarcs}\label{app:reverse}

In this appendix, we extend the multiperiod hypergraph construction framework to account for reversible hyperarcs. Let $r(\mathcal{A}) \subseteq \mathcal{A}$ denote the set of reversible hyperarcs:
$$
r(\mathcal{A})= \{a \in \mathcal{A}: \exists a' \in \A \mbox{ with } S_a=T_{a'} \text{ and } S_{a'}=T_a\}.
$$
For each $a \in r(\mathcal{A})$, we denote by $r(a)$ its reverse hyperarc, defined by exchanging its source and target sets, i.e., $\mathbb{S}_{\cdot a} = \mathbb{T}_{\cdot r(a)}$ and $\mathbb{T}_{\cdot a} = \mathbb{S}_{\cdot r(a)}$.

Under the general framework where $r(\mathcal{A}) \neq \emptyset$, the synergistic flow associated with a hyperarc $a \in \mathcal{A}$ is defined as
\begin{align}
f_a =
\begin{cases}
\kappa_a \displaystyle\prod_{v \in \mathcal{N}} x_v^{\mathbb{S}_{va}}, 
& \forall a \in \mathcal{A} \setminus r(\mathcal{A}),\\[6pt]
\kappa_a \displaystyle\prod_{v \in \mathcal{N}} x_v^{\mathbb{S}_{va}} 
- \kappa_{r(a)} \displaystyle\prod_{v \in \mathcal{N}} x_v^{\mathbb{T}_{va}}, 
& \forall a \in r(\mathcal{A}),
\end{cases}
\label{ctr:12_0_app}
\end{align}
which captures the net effect of forward and reverse contributions. This formulation is consistent with generalized kinetic representations in reaction systems where reverse transformations act as competing processes~\citep[see, e.g.][]{results_anderson2019}.

Introducing the binary activation variables $z_a^t$ and state variables $x_v^t$, the corresponding multiperiod representation becomes
\begin{align}
f_a^t =
\begin{cases}
\kappa_a z_a^t \displaystyle\prod_{v \in \mathcal{N}} (x_v^t)^{\mathbb{S}_{va}}, 
& \forall a \in \mathcal{A} \setminus r(\mathcal{A}), \; t \in \mathcal{T},\\[6pt]
\kappa_a z_a^t \displaystyle\prod_{v \in \mathcal{N}} (x_v^t)^{\mathbb{S}_{va}} 
- \kappa_{r(a)} z_a^t \displaystyle\prod_{v \in \mathcal{N}} (x_v^t)^{\mathbb{T}_{va}}, 
& \forall a \in r(\mathcal{A}), \; t \in \mathcal{T}.
\end{cases}
\label{ctr:12_app}
\end{align}

\subsection{Logarithmic reformulation}

To derive a tractable mixed-integer formulation, we introduce auxiliary variables for the reverse contribution of reversible hyperarcs. Let $\tilde{f}_a^t$ denote the reverse flow associated with $a \in r(\mathcal{A})$. Then, the following constraints model the reverse component:

\begin{align*}
\tilde g_{a}^t &\geq \log(\kappa_{r(a)}) + \sum_{v\in \mathcal{N}} \mathbb{T}_{va} h_v^t - \Delta_a (1-z_a^t), 
\quad \forall a \in r(\mathcal{A}), \; t \in \mathcal{T},\\
\tilde g_{a}^t &\leq \log(\kappa_{r(a)}) + \sum_{v\in \mathcal{N}} \mathbb{T}_{va} h_v^t + \Delta_a (1-z_a^t), 
\quad \forall a \in r(\mathcal{A}), \; t \in \mathcal{T},\\
\tilde g_{a}^t &= \log(\tilde f_{a}^t - z_a^t + 1), 
\quad \forall a \in r(\mathcal{A}), \; t \in \mathcal{T},\\
\tilde f_{a}^t &\leq \Delta_a z_a^t, 
\quad \forall a \in r(\mathcal{A}), \; t\in \mathcal{T},\\
\tilde f_{a}^t &\geq \varepsilon_a z_a^t, 
\quad \forall a \in r(\mathcal{A}), \; t\in \mathcal{T}.
\end{align*}

This construction isolates the reverse contribution, allowing the net flow to be represented as the difference between forward and reverse components, while preserving the logarithmic linearization scheme used in the main model.

The presence of reverse flows modifies both the feasibility and dynamic constraints. In particular, constraints \eqref{ctr:7} and \eqref{ctr:9} are replaced by

\begin{align*}
\sum_{a \in \mathcal{A}} \mathbb{Q}_{va} f_a^t 
- \sum_{a \in r(\mathcal{A})} \mathbb{Q}_{va} \tilde f_a^t 
&\geq \varepsilon - \Delta(1-y_v^t), 
\quad \forall v \in \mathcal{N}, \; t\in \mathcal{T},\\
x_v^t &= x_v^{t-1} 
+ \sum_{a \in \mathcal{A}} \mathbb{Q}_{va} f_a^t 
- \sum_{a \in r(\mathcal{A})} \mathbb{Q}_{va} \tilde f_a^t, 
\quad \forall v \in \mathcal{N}, \; t \in \mathcal{T}.
\end{align*}

These expressions ensure that the net production accounts for both forward and reverse contributions, preserving mass balance and consistency of the multiperiod dynamics.

This extension preserves the core structure of the model while allowing for richer interaction patterns. From a modeling perspective, reversible hyperarcs introduce competing mechanisms that may attenuate or reinforce amplification processes, a phenomenon closely related to autocatalytic and inhibitory interactions in reaction networks~\citep{peng_hierarchical_2022}. Importantly, the separation of forward and reverse components maintains tractability and enables the direct incorporation of reversibility into the optimization framework.

\end{document}